\documentclass[11pt]{amsart}

\usepackage{amssymb}
\usepackage[all]{xy}
\usepackage[colorlinks=true, urlcolor=rltblue, citecolor=drkgreen, linkcolor=drkred] {hyperref}
\usepackage{color}
\usepackage{pdfsync}
\definecolor{rltblue}{rgb}{0,0,0.4}
\definecolor{drkred}{rgb}{0.6,0,0}
\definecolor{drkgreen}{rgb}{0,0.4,0}
\usepackage{graphicx}
\newtheorem{thm}{Theorem}[section]
\newtheorem{lemma}[thm]{Lemma}

\newtheorem{theorem}[thm]{Theorem}
\newtheorem{corollary}[thm]{Corollary}

\theoremstyle{definition}
\newtheorem{definition}[thm]{Definition}

\theoremstyle{remark}

\newtheorem{obs}[thm]{Observation}
\newtheorem{remark}[thm]{Remark}

\newtheorem{historic}[thm]{Historic Remark}

\theoremstyle{plain}

\newcounter{contenumi}

\def\teq{\equiv_T}

\def\upto{\mathop{\upharpoonright}}
\def\la{\langle}
\def\ra{\rangle}
\def\and{\mathrel{\&}}

\def\sminus{\smallsetminus}

\def\isom{\cong}
\def\Si{\Sigma}
\newcommand\rightdate[1]{\footnotetext{  Saved: #1 \\ Compiled: \today}}
\def\A{\mathcal{A}}
\def\B{\mathcal{B}}
\def\C{\mathcal{C}}

\def\om{\omega}
\def\bbar{\bar{b}}

\def\si{\sigma}
\def\b{\beta}
\def\a{\alpha}
\def\b{\beta}

\newcommand \converges{\downarrow}

\def\A{\mathcal A}

\def\cbar{\bar{c}}

\def\abar{\bar{a}}
\def\bbar{\bar{b}}

\def\dbar{\bar{d}}

\def\itt{{\mathtt {in}}}

\def\Pii{\Pi^\itt}

\def\concat{^\frown}
\def\g{\gamma}
\def\d{\delta}

\def\Si{\Sigma}

\def\om{\omega}

\def\implies{\Rightarrow}

\newcommand{\J}{\mathcal{J}}

\newcommand{\fseq}[1]{2^{\circ #1}}
\newcommand{\Bai}{\ensuremath{\om^\om}}
\newcommand{\Seq}{{\ensuremath{\om^{<\om}}}}
\newcommand{\diverges}{\uparrow}
\newcommand{\Jom}[1]{J^{\om^{#1}}}

\def\Nleq{\trianglelefteq}

\def\Aleq{\preccurlyeq}

\def\dmi{\ \ {}^\cdot\!\!\!\!\! -\ }
\def\basi{\bar{\si}}
\def\batau{\bar{\tau}}
\title{Priority arguments via true strages}

\author{Antonio Montalb\'an}
\thanks{The author was partially supported by NSF grant DMS-0901169 and the Packard Fellowship.
The author would also like to thank Mariya Soskova and Uri Andrews for proof reading this paper.
This paper was written while the author participated in the Buenos Aires 
Semester in Computability, Complexity and Randomness, 2013.
}

\address{Department of Mathematics\\
University of California, Berkeley}

\email{antonio@math.berkeley.edu}
\urladdr{\href{http://www.math.berkeley.edu/~antonio/index.html}{www.math.berkeley.edu/$\sim$antonio}}

\begin{document}

\rightdate{June 11, 2013 -- submitted}
\maketitle

%t\tableofcontents

%

\begin{abstract}
We describe a variation of Ash's $\eta$-system, and give a new proof of Ash's metatheorem.
As an application, we prove a generalization of Ash and Knight's theorem on pairs of structures.
\end{abstract}

%%%%%%%%%%%%%%%%%%%%%%%%%%%%%%%%%%%%%%%%%%%%%%%%%%%%%%%%%%%%%%%%%%%%%%%%%%%%%%%%%%%%%%%%%%%%%%%%%%%%%%%%%%%%%%%%
\section{Introduction}

Ash's metatheorem is one of the gems of computability theory.
It gives a general framework for doing $0^{(\eta)}$-priority constructions which can be applied to a whole variety of cases, mostly coming from Computable Structure Theory. 
What makes it so special is that with a relatively simple combinatorial condition it captures the essence of $0^{(\eta)}$-priority constructions.
It was first introduced by Ash in 1986 \cite{Ash86hyp, Ash86}, and several slightly different versions where proposed later in the 90's by Ash and Knight \cite{Ash90, AK94Ramified, AK94Mixed,Kni95}.
The standard reference for Ash's metatheorem and its applications is Ash and Knight's book \cite{AK00}.

The objective of this paper is twofolded.
First, we present a new version of Ash's metatheorem and give a new proof for it.
Our new definition of an $\eta$-system is in terms of the stage-by-stage construction that it produces, which we hope might help the readers who are used to thinking in these terms.
The key notion in this new presentation is that of a $\xi$-true stage; these are the stages where our approximations to the $\xi$'th jump are correct.
Most computability theorist are used to proofs that use true stages for the enumeration of $0'$, and they might find the extensions to higher ordinals quite natural. 
The proof is not really simpler, but the difficult part of the proof is now separated from the metatheorem, and moved inside the definitions of the approximation of the jumps, and of ``apparent $\xi$-true stages.''
In Section \ref{se: true stages} we develop a way of producing, at each finite stage, an approximation to the $\xi$-jump.
The ideas to define  such approximations are not completely new, and different approaches can be found in the literature. 
The approximations of the jumps we use are based on definitions from Marcone and Montalb\'an \cite{MMVeblen}, because of their nice combinatorial properties that will be useful to us, although some important modifications are needed.
Assuming some basic properties for these notions, the proofs of the metatheorems in Sections \ref{se: n systems} and \ref{se: eta systems} are quite simple, and they concentrate on the combinatorial aspects of the $\eta$-system, which are now clearly separated from the combinatorial aspects of the approximations to the jumps.

Depending on the application at hand, one framework might me easier to apply than another.
Let us remark that there are two other known general frameworks for priority arguments:
The {\em workers method} developed by Harrington \cite{Har76,Har79} (see also Knight \cite{Kni90a,Kni90b}), and the {\em iterated trees of strategies method}, developed by Lempp and Lerman \cite{LL90,LL95,Ler10}.

Second, we give a particular application of our new version of Ash's metatheorem.
Neither a standard Ash's $\eta$-system, nor our $\eta$-system defined in Section \ref{se: eta systems} seems to be enough for this application, and we thus need to provide a more general version in Section \ref{se: small modification}.
This modification is rather small, but it seems to give the $\eta$-system much extra power.
We remark that this modification seems small only because of the presentation we chose, and it is not clear how to state this modification using Ash's original framework.

Our application, Theorem \ref{thm: app}, is a generalization of Ash and Knight's theorem on pairs of structures \cite{AK90} which says the following:
Suppose we have a $\Si^0_\eta$ question, and we want to build, uniformly, a computable structure $\C$ which is isomorphic to a given structure $\A_1$ if the answer to the question is `yes', and isomorphic to another given structure $\A_0$ otherwise.
Then, we can do this if we know that  $\A_0\geq_\eta\A_1$ (where $\leq_\eta$ is the $\eta$-th back-and-forth relation) and that $\A_0,\A_1$ satisfy an effectiveness condition called $\eta$-friendliness.
This result has been used in all sorts of applications around the literature.
An important corollary is that $\A_0\geq_\eta\A_1$ if and only if the problem of distinguishing between $\A_0$ and $\A_1$ (given a structure that we know is isomorphic to one of the two) is boldface ${\mathbf \Si^0_\eta}$-hard (see \cite[Theorem 2.2]{HMapprox}). 

In our generalization, instead of having two structures we have a binary tree of structures.
This tree has ordinal height, and if we branch at the $\xi$-th level of the tree, then the structures on one side of the branch are $\geq_\xi$-greater than the ones on the other branch.
The theorem then says that if we are given a branch, where calculating the $\xi$-th bit of the branch is $\Si^0_\xi$, then we can uniformly build a computable structure isomorphic to the one lying at that branch of the tree.
We also assume the branch has finitely many 1's, and again, we assume the structures in the tree are all $\eta$-friendly.
See Section \ref{se: app} for the full statement and definitions. 

We are  interested in this particular application because it is a key lemma in an upcoming paper \cite{MonIntermediateVC} where the author gives a Vaught's-conjecture type characterization of the classes of structures which are intermediate for computable embeddability, which is the effective version of the Friedman-Stanley reducibility \cite{FS89} on classes of structures studied on \cite{FF09, FFHKMM}.

%%%%%%%%%%%%%%%%%%%%%%%%%%%%%%%%%%%%%%%%%%%
\section{$n$-systems} \label{se: n systems}

We start by defining $n$-systems for finite $n$, just because it is easier to start by picturing the finite case, and because some readers might only be interested in the finite case.

Our  $n$-systems will rely on the notion of ``apparent $n$-true stage.''
We shall define this notion precisely in Section \ref{see: xi-true stages}.
For now we just list the properties this notion should have, and see how these properties are enough to prove our metatheorem. 
Since most computability theorist already have a good intuition as to what this notion should be, we expect the reader to be able to read this section without having a concrete definition in mind.
Otherwise, the reader can skip to Sections \ref{se: true stages} and  \ref{se: belief}, and only then come back to this section.

%%%%%%%%%%%%%%%%
\subsection{Apparent true stages}  \label{see: apparent true}

As we know, a true stage is one at which we are guessing a finite initial segment of $0'$ correctly.
We will call them {\em 1-true stages}.
Similarly, an $n$-true stage is one at which we are guessing a finite initial segment of $0^{(n)}$ correctly.
%A stage $s$ {\em looks $n$-true at} a later stage $t$, if the guess for $0^{(n)}$ at $s$ is compatible with the guess at $t$.

Let us define the notation we will use.
For each computable ordinal $\xi$ we will use $\nabla^\xi\in\om^\om$ to denote a $\Delta^0_\xi$-complete function (Definition \ref{def: nabla}).
Notice that for finite $n$, $\nabla^{n+1}$ is Turing equivalent to $0^{(n)}$, and for infinite $\xi$, $\nabla^\xi\teq 0^{(\xi)}$.
(It might help the reader to start by considering only the case when $\xi$ is just a finite number.)
For each $\xi$, and at each stage $s$, we will computably define a finite string $\nabla^\xi_s$ which is trying to guess an initial segment of $\nabla^\xi$.
These approximations will satisfy the following properties, that we will prove in Lemmas \ref{le: seq of xi true}, \ref{le: finitely many} and \ref{le: nabla nested}.

\begin{enumerate}\renewcommand{\theenumi}{N\arabic{enumi}} 
\item For each $\xi$, there is a sub-sequence $\{t_i:i\in\om\}$ such that  $\bigcup_{i\in\om}\nabla^\xi_{t_i} =\nabla^\xi$.  \label{pa:N1}
\item For every $s$, there are only finitely many $\xi$'s with $\nabla^\xi_s\neq\la\ra$ (where $\la\ra$ is the empty sequence).  \label{pa:N2}
\item For $\g\leq \xi$ and $s\leq t$, if $\la\ra\neq\nabla^\xi_s\subseteq\nabla^\xi_t$, then $\nabla^\g_s\subseteq\nabla^\g_t$.  \label{pa:N3}
\end{enumerate}

In Section \ref{see: xi-true stages} we will define a sequence of relations $(\leq_\xi)_{\xi\leq\eta}$ on $\om$ that will represent the notion of apparent $\xi$-true stage.
That is, for $s<t\in\om$, we say that  {\em $s$ looks $\xi$-true at $t$}, or that $t$ {\em $\xi$-believes in} $s$, if $s\leq_\xi t$.
For finite $n$ we can let $s\leq_n t$ if and only if $\nabla^{n+1}_s\subseteq\nabla^{n+1}_t$, and we would get all the properties we need.
However, for transfinite $\xi$, the relations $\nabla^{\xi+1}_s\subseteq\nabla^{\xi+1}_t$ fail to have some of these properties, and we will have to define $\leq_\xi$ slightly differently.
For now, all we need to know is that the relations $\leq_\xi$ on $\om$ will satisfy the following  properties (Lemma \ref{le: xi-belief final}):

\begin{enumerate}\renewcommand{\theenumi}{B\arabic{enumi}} \setcounter{enumi}{-1}
\item $\leq_0$ is just the standard ordering on $\om$.   \label{pa: 0}
\item The relations $\leq_\xi$ are uniformly computable. \label{pa: 0.5}
\item Each $\leq_\xi$  is a {\em pre-ordering} (i.e., reflexive and transitive).   \label{pa: 1}
\item The sequence of relations is {\em nested} (i.e., if $\g\leq \xi$ and $s\leq_\xi t$, then $s\leq_\g t$).\label{pa: 2}
\item The sequence of relations is {\em continuous} (i.e., if $\lambda$ is a limit ordinal, then $\leq_\lambda=\bigcap_{\xi<\lambda}\leq_\xi$).\label{pa:4}
\item For every $s<t\in\om$, $s\leq_\xi t \implies  \nabla^{\xi+1}_s\subseteq\nabla^{\xi+1}_t$.    \label{pa:5}
\item For every $\xi$, there exist stages $t_0\leq_\xi t_1\leq_\xi\cdots$ with  $\bigcup_{i\in\om}\nabla^{\xi+1}_{t_i} =\nabla^{\xi+1}$. \label{pa:6}
%\item    \label{pa: 4}
%given $r< s< t$, if $r\leq_{n+1} t$ and $s\leq_n t$, then $r\leq_{n+1} s$.  
\end{enumerate}

But the most important property the relations $\leq_\xi$ satisfy is the following:

\begin{enumerate}\renewcommand{\theenumi}{$\clubsuit$} 
\item For every $\xi$, and every $r< s< t$, if $r\leq_{\xi+1} t$ and $s\leq_\xi t$, then $r\leq_{\xi+1} s$.
\end{enumerate}

To see why $(\clubsuit)$ should be true, suppose that  $r<s$ both have compatible guesses for $\nabla^{\xi+1}$, but $r$ does not look $(\xi+1)$-true at $s$.
The only way this could happen is that $s$ witnessed some convergence with oracle $\nabla^{\xi+1}$ which $r$ had not seen yet.
This convergence stays present thereafter, so if a later stage $t$ has the same guesses for $\nabla^{\xi+1}$, it will see the same convergence and we will have $r\not\leq_{\xi+1} t$.

\begin{obs}\label{rm: tree like}
Given $r<s<t$, if $r\leq_\xi t$ and $s\leq_\xi t$, then $r\leq_\xi s$:
For successor $\xi$ use $(\clubsuit)$ and that $s\leq_{\xi-1} t$.
For limit $\xi$, we have that for every $\g<\xi$, we can use $(\clubsuit)$, and that $r\leq_{\g+1} t$ and $s\leq_\g t$, to get $r\leq_{\g+1} s$.
By continuity we get that $r\leq_\xi s$.
\end{obs}

%%%%%%%%%%%%%%%%%%%
\subsection{The system}

We are now ready to define $n$-stystems.

\begin{definition}
An $n$-system consists of a triple $(L, P_0, (\leq^L_i)_{i\leq n})$ where
\begin{itemize}
\item $L$ is a computable subset of $\om$, called the set of {\em states}.
\item $P_0$ is a computable subtree of $L^{<\om}$, called the {\em action tree}.
\item $(\leq^L_i)_{i\leq n}$ is a computable nested sequence of pre-orderings on $L$, called the {\em restraint relations}.
\end{itemize}
\end{definition}

\begin{remark}
For the reader familiar with Ash's $n$-systems, let us note that this is what he would have called an $(n+1)$-system, rather than an $n$-system.
The alternating tree $P$ in Ash's system corresponds to the restriction of  $P_0$ to the $n$-true stages.
Instead of the instruction function $q$ in Ash's system, we will use the approximations to the jumps and code it directly within $P_0$.
The enumeration function $E$ is defined below, and plays the same role as in Ash's systems.
\end{remark}

At each stage $s$ of the construction we will choose $\ell_{s}\in L$ which determines the state of the construction so far.
We must always have $\la \ell_0,...,\ell_s\ra \in P_0$.
If we were to do an analogy with forcing constructions, the elements of $L$ would be called ``forcing conditions.''
The action tree must guarantee that the state $\ell_{s}$ built at stage $s$ is compatible with $\nabla^{n+1}_s$, our current guess for $0^{n}$.
%Let us remark that, using the recursion theorem, we can assume that $\nabla^{n}$ not only knows about the jumps of $0$, but also about the jumps of the construction.
%So, using the current guesses for the jumps of the construction, the action tree $P_0$ will enforce some action is taken to satisfy the requirements for the construction.

The restraint relations, $\leq^L_i$, are the ones that guarantee that not much injury is done at stages where our guesses for the jumps are wrong.
The main issue in $0^{(n)}$-priority constructions is that at a given stage one might  guess wrongly which requirements need to be respected and which do not, and this is why the injuries need to be carefully controlled in a way that we can recover from them later. 
The restraints are going to come in different levels.
Ideally, an $i$-restraint would like all the upcoming states in the construction to be $\geq_i^L$-greater than the current one.
But this is too much to ask in general, and we will only require it on the $i$-true stages.
To get this, the attitude of the construction is the following: If $s$ $i$-believes in a previous stage $t$ (i.e., $s\geq_i t$), then $s$ must respect the $i$-restraint imposed by $t$ (i.e., $\ell_{s}\geq^L_i \ell_t$).

The objective of many constructions in computability theory is to build a c.e.\ object of some kind, either a collection of c.e.\ sets, or a diagram of a structure, or a Turing operator, etc.
Each $n$-system comes together with an {\em enumeration function} in charge of enumerating the outcome of the construction.
An enumeration function is  a c.e.\ set $E\subseteq L\times \om$. 
We denote $E(l)=\{k\in\om: (l,k)\in E\}$. 
It must satisfy that for $\ell_0,\ell_1\in L$, if $\ell_0\leq^L_0\ell_1$, then $E(\ell_0)\subseteq E(\ell_1)$.
The idea is that when we play $\ell$ at some stage in our construction, we enumerate $E(\ell)$ in our outcome set.

The objective of the construction is to define a computable infinite sequence $\pi=(\ell_0,\ell_1,....)$ with the following property. 

\begin{definition}
A {\em 0-run} for an $n$-system $(L, P_0, (\leq^L_i)_{i\leq n})$ is a finite or infinite sequence $\pi=(\ell_0,\ell_1,....)\in P_0\cup[P_0]$ (where $[P_0]$ is the set of paths through $P_0$) such that, for each $s,t<|\pi|$ and $i\leq n$, we have that
\[
s \leq_i t  \ \  \ \implies \ \ \  \ell_{s} \leq^L_i \ell_{t}.
\]
Given a 0-run, we let $E(\pi)=\bigcup_{i<|\pi|}E(\ell_i)$.
\end{definition}

The 0-run $\pi$ is what we call ``the construction.''
Why do we desire to build an infinite computable 0-run $\pi$?
Let $t_0,t_1,....$ be the sub-sequence of $n$-true stages (i.e.\ the stages $t$ as in (\ref{pa:6}), for which $\nabla^{n+1}_t$ is correct).
We call $\pi_{n}=\la \ell_{t_0},\ell_{t_1},...\ra$  an {\em $n$-run}.
(This is closer in essence to what Ash calls a run.)
We then have that 
\[
\ell_{t_0}\leq^L_{n} \ell_{t_1}\leq^L_{n} \ell_{t_2}\leq^L_{n} ....
\]
and furthermore, each $\ell_{t_i}$ was buit under the correct assumption about the construction.
So any requirement that is satisfied along this sequence, should be correctly satisfied forever.
On the other hand, note that 
\[
E(\pi) = E(\pi_{n}) = \bigcup_{i\in\om}E(\ell_{t_i}).
\]
Even though $E(\pi)$ is just c.e., it can be defined in terms of what is enumerated along the $n$-true stages under the correct guess for $\nabla^{n+1}$.
This should guarantee that  $E(\pi)$ satisfies the desired requirements.

%%%%%%%%%%%%%%%%%%%%%
\subsection{The Metatheorem}

The theorem  below gives a sufficient condition  to build a computable infinite $0$-run on an $n$-system.
This is the condition that controls the injuries of the construction, the one that allows us to recover from the injuries done at stages where we had incorrect information.

\begin{definition}
We say that an $n$-system $(L, P_0, (\leq^L_i)_{i\leq n})$ satisfies the {\em extendibility condition} if:
whenever we have a finite 0-run $p=\la \ell_0,...,\ell_{s-1}\ra$, and a finite sequence of $s_k \leq s_{k-1} \leq ....\leq s_0< s$  with $k\leq n$ such that, for all $i$, $\ell_{s_i} \leq^L_i \ell_{s_{i-1}}$, there exists an $\ell\in L$ such that $p\concat\ell\in P_0$ and $\ell_{s_i} \leq^L_i \ell$ for all $i\leq k$. (See diagram below.)
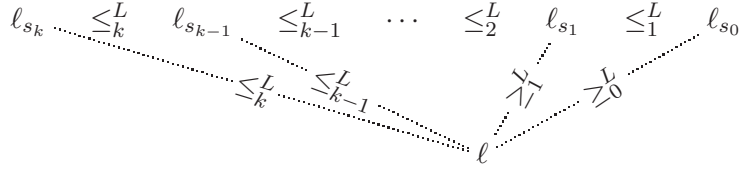
\[
\xymatrix@C=10pt@R=7pt{
%	& 			& 		& 			&		& 	 
%	s \ar@{.}|{\rotatebox[origin=c]{15}{$\leq_{k}$}}[lllllddd] \ar@{.}|{\rotatebox[origin=c]{22}{$\leq_{k-1}$}}[lllddd]  \ar@{.}|{\rotatebox[origin=c]{-60}{$\geq_{{1}}$}}[rddd]  \ar@{.}|{\rotatebox[origin=c]{-27}{$\geq_{0}$}}[rrrddd]  \\ \\ \\ 
%{{s_k}}  &  \leq_{k} &   {s_{k-1}} 	& \leq_{{k-1}} &  \cdots &  \leq_{2}  &  {{s_1}}  &  \leq_{1}    &  {{s_0}}  \\ 	
{\ell_{s_k}}  & \leq^L_{k} &   \ell_{s_{k-1}} 	& \leq^L_{{k-1}} &  \cdots &  \leq^L_{2}   &  {\ell_{s_1}}  & \leq^L_{1}    &  {\ell_{s_0}}  \\ 	\\
\\
	& 			& 		& 			&		& 	  
	{\ell} \ar@{.}|{\rotatebox[origin=c]{-17}{$\leq^L_{k}$}}[llllluuu] \ar@{.}|{\rotatebox[origin=c]{-20}{$\leq^L_{k-1}$}}[llluuu]  \ar@{.}|{\rotatebox[origin=c]{62}{$\geq^L_{{1}}$}}[ruuu]  \ar@{.}|{\rotatebox[origin=c]{34}{$\geq^L_{0}$}}[rrruuu] 
}\]
%\[
%\xymatrix@C=10pt@R=7pt{
%%	& 			& 		& 			&		& 	 
%%	s \ar@{.}|{\rotatebox[origin=c]{15}{$\leq_{k}$}}[lllllddd] \ar@{.}|{\rotatebox[origin=c]{22}{$\leq_{k-1}$}}[lllddd]  \ar@{.}|{\rotatebox[origin=c]{-60}{$\geq_{{1}}$}}[rddd]  \ar@{.}|{\rotatebox[origin=c]{-27}{$\geq_{0}$}}[rrrddd]  \\ \\ \\ 
%%{{s_k}}  &  \leq_{k} &   {s_{k-1}} 	& \leq_{{k-1}} &  \cdots &  \leq_{2}  &  {{s_1}}  &  \leq_{1}    &  {{s_0}}  \\ 	
%{\ell_{k}}  & \leq^L_{k} &   \ell_{{k-1}} 	& \leq^L_{{k-1}} &  \cdots &  \leq^L_{2}   &  {\ell_{1}}  & \leq^L_{1}    &  {\ell_{0}}  \\ 	\\
%\\
%	& 			& 		& 			&		& 	  
%	{\ell} \ar@{.}|{\rotatebox[origin=c]{-17}{$\leq^L_{k}$}}[llllluuu] \ar@{.}|{\rotatebox[origin=c]{-20}{$\leq^L_{k-1}$}}[llluuu]  \ar@{.}|{\rotatebox[origin=c]{62}{$\geq^L_{{1}}$}}[ruuu]  \ar@{.}|{\rotatebox[origin=c]{34}{$\geq^L_{0}$}}[rrruuu] 
%}\]
\end{definition}

%(As is drawn in the picture, we notice that $s_i\leq_is_{i-1}$.
%This follows from property $(\clubsuit)$ using that  $s_i\leq_{i}s$ and $s_{i-1}\leq_{i-1}s$.)

The following is the main lemma about $n$-systems.
Ash calls this the ``metatheorem.''
Its proof is quite simple, due to the carefully crafted definitions above.
What makes it the main lemma is that it shows how the extendibility condition is the combinatorial core of the construction.

\begin{theorem}
For every $n$-system $(L, P_0, (\leq^L_i)_{i\leq n})$ with the extendibility condition, there exists an infinite computable  0-run $\pi$.
Furthermore, a $0$-run can be built uniformly in the $n$-system.
\end{theorem}
\begin{proof}
We build $\pi$ stage by stage.
To start, the extendibility condition trivialy gives us some $\ell_0$ such that $\la\ell_0\ra\in P_0$, which is enough to have a 0-run of length one.
Suppose we have already built a finite 0-run $p=\la \ell_0,...,\ell_{s-1}\ra$ and we want to define $\ell_{s}$.

We need to find $\ell\in L$ such that $p\concat\ell\in P_0$ and such that, for every $t< s$ and $i\leq n$, if $t\leq_i s$, then $\ell_t\leq^L_i \ell_{s}$.
To find such an $\ell$ it is not necessary to consider all such $t$'s, but only the maximal among them for each $i$.
For $i\leq n$, let $s_i$ be the largest $t<s$ such that $t\leq_i s$, if such a $t$ exists.
Let $k\leq n$ be the largest such that $s_k$ exists.
Note that for $i<j\leq n$, if $s_j$ exists then so does $s_i$ because  $s_j\leq_i s$ (by the nesting of the relations $\leq_i$).
Furthermore, we have that $s_{j}\leq s_{i}$.
Actually, using $(\clubsuit)$ and that   $s_i\leq_{i}s$ and $s_{i-1}\leq_{i-1}s$, we have that $ s_i\leq_i s_{i-1}$.
And hence, since $\pi$ is a 0-run, we have that $\ell_{s_i}\leq^L_i\ell_{s_{i-1}}$ for all $i$.
The extendibility condition then gives us an $\ell\in L$ such that $p\concat\ell\in P_0$ and $\ell_{s_i} \leq^L_i \ell$ for all $i$.
Let $\ell_{s}=\ell$.
We now claim that $\la \ell_0,...,\ell_{s-1},\ell_{s}\ra$ is a 0-run:
Suppose that $t< s$ is such that $t\leq_i s$.
By the definition of $s_i$ we get that $t\leq s_i$, and by Observation \ref{rm: tree like} we get that $t\leq_is_i$, and hence that $\ell_t\leq^L_i \ell_{s_i}$.
By the transitivity of $\leq^L_i$ we get that $\ell_t\leq^L_i \ell_{s}$.
\end{proof}

%%%%%%%%%%%%%%%%%%%%%%%%%%%%%%%
\section{$\eta$-systems}   \label{se: eta systems}

The generalization from finite $n$ to transfinite $\eta$ is not too drastic.
The definitions of $\eta$-system and of 0-run are exactly the same as above, just thinking of $n$ as a transfinite ordinal rather than a natural number.
The extendibility condition needs to be slightly modified.

Assume  that $\eta$ is a computable well-ordering, and that its presentation is nice enough so that the successor function is computable.
When we write $\xi\leq\eta$, we mean that $\xi$ is a member of this particular presentation of $\eta+1$.

%For the relations of ``$i$-believes in'' we will require an extra condition at the limit levels.
%We define this pre-orders in Section \ref{see: xi-true stages}.
%For now we assume $(\leq_\xi)_{\xi\leq\eta}$ is a computable nested sequence of pre-orders satisfying property $(\clubsuit)$.
%All these properties are defined exactly as in the previous section just thinking of $n$ as a transfinite ordinal.
%Not only we require $(\leq_\xi)_{\xi\leq\eta}$ to be a nested sequence, but also we require it is {\em continuously nested}:
%That is, for limit $\lambda\leq\eta$, $\leq_\lambda=\bigcap_{\xi<\lambda}\leq_\xi$.
%
%We also assume that $\leq_0=\leq$, and that $s\leq_\eta t\iff\nabla^{\eta+1}_s\subseteq\nabla^{\eta+1}_t$, provided $\la\ra\neq\nabla^{\eta+1}_s$.
%This implies that there is an infinite sequence of $\eta$-true stages $t_1\leq_\eta t_2\leq_\eta\cdots$, on which the guesses to $\nabla^{\eta+1}$ are all correct.

\begin{definition}
We say that an $\eta$-system $(L, P_0, (\leq^L_\xi)_{\xi\leq\eta})$ satisfies the {\em extendibility condition} if:
whenever we have a finite 0-run $p=\la \ell_0,...,\ell_{s-1}\ra$, stages $s_k < s_{k-1} < ....< s_0< s$, and ordinals $\xi_0<\xi_1<\cdots<\xi_k\leq\eta$ such that, for all $i< k$,  $\ell_{s_{i+1}} \leq^L_{\xi_i+1} \ell_{s_i}$, there exists an $\ell\in L$ such that $p\concat\ell\in P_0$ and, for all $i\leq k$, $\ell_{s_i} \leq^L_{\xi_i} \ell$.
(See diagram below.)
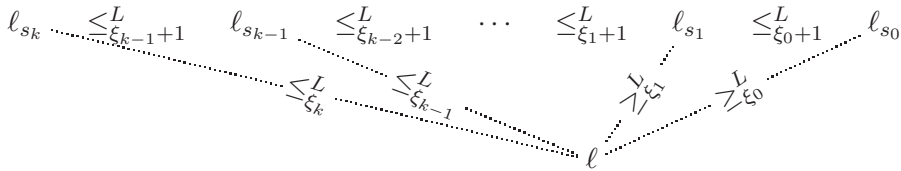
\[
\xymatrix@C=10pt@R=7pt{
%	& 			& 		& 			&		& 	 
%	s \ar@{.}|{\rotatebox[origin=c]{14}{$\leq_{\xi_k}$}}[lllllddd] \ar@{.}|{\rotatebox[origin=c]{21}{$\leq_{\xi_{k-1}}$}}[lllddd]  \ar@{.}|{\rotatebox[origin=c]{-54}{$\geq_{\xi_{1}}$}}[rddd]  \ar@{.}|{\rotatebox[origin=c]{-24}{$\geq_{\xi_0}$}}[rrrddd]  \\ \\ \\ 
%{{s_k}}  &  \leq_{\xi_{k-1}+1} &   {s_{k-1}} 	& \leq_{\xi_{k-2}+1} &  \cdots &  \leq_{\xi_1+1}  &  {{s_1}}  &  \leq_{\xi_0+1}    &  {{s_0}}  \\ 	
{\ell_{s_k}}  & \leq^L_{\xi_{k-1}+1} &   \ell_{s_{k-1}} 	& \leq^L_{\xi_{k-2}+1} &  \cdots &  \leq^L_{\xi_1+1}   &  {\ell_{s_1}}  & \leq^L_{\xi_0+1}    &  {\ell_{s_0}}  \\ 	\\
\\
	& 			& 		& 			&		& 	  
	{\ell} \ar@{.}|{\rotatebox[origin=c]{-17}{$\leq^L_{\xi_k}$}}[llllluuu] \ar@{.}|{\rotatebox[origin=c]{-20}{$\leq^L_{\xi_{k-1}}$}}[llluuu]  \ar@{.}|{\rotatebox[origin=c]{58}{$\geq^L_{\xi_{1}}$}}[ruuu]  \ar@{.}|{\rotatebox[origin=c]{34}{$\geq^L_{\xi_0}$}}[rrruuu] 
}\]
\end{definition}

%As is drawn in the picture, we notice that $s_{i+1}\leq_{\xi_{i}+1}s_{i}$.
%This follows from property $(\clubsuit)$ using that $s_{i}\leq_{\xi_{i}}s$, and that $s_{i+1}\leq_{\xi_{i+1}}s$ which implies that  $s_{i+1}\leq_{\xi_{i}+1}s$.

\begin{theorem}\label{thm: meta eta}
For every $\eta$-system $(L, P_0, (\leq^L_\xi)_{\xi\leq \eta})$ with the extendibility condition, there exists a computable infinite 0-run $\pi$.
Furthermore, a $0$-run can be build uniformly in the $n$-system.
\end{theorem}
\begin{proof}
Let $\ell_0$ be such that $\la\ell_0\ra\in P_0$, which exists by the trivial case in the extendibility condition.
Suppose we have already built a finite 0-run $p=\la \ell_0,...,\ell_{s-1}\ra$ and we want to define $\ell_s$.

We want to define $\ell_{s}\in L$ such that $p\concat \ell_s\in P_0$ and, for every $t<s$ and every $\xi\leq\eta$, if $t\leq_\xi s$ then $\ell_t\leq^L_{\xi}s$.
Let $\eta_0\leq\eta$ be the largest such that there exists some $t<s$ with $t\leq_{\eta_0} s$.
A largest such ordinal exists by the continuity of $(\leq_\xi)_{\xi\leq\eta}$.
For each $\xi\leq\eta_0$, let $t_\xi< s$ be the largest $t$ such that $t\leq_\xi s$.
Notice that if $t\leq_\xi s$, then $t\leq t_\xi$ and by Observation \ref{rm: tree like}, $t\leq_\xi  t_\xi\leq_\xi s$.
So it is enough to get $\ell$ with $p\concat\ell\in P_0$ such that, for every $\xi\leq\eta_0$, $\ell_{t_\xi}\leq^L_{\xi}\ell$.

There are infinitely many $\xi$'s, but only finitely many possible values for $t_\xi< s$, so they  must repeat a lot.
Using that the relations $(\leq_\xi)_{\xi\leq\eta}$ are nested, we see that if $\xi\leq\zeta\leq\eta_0$ then $t_\zeta\leq t_\xi$.
We now want to define stages $s_k<....<s_0<s$ so that $\{s_0,...,s_k\}=\{t_{\xi}:\xi\leq\eta_0\}$ as sets, but we need to define them effectively.
Let $s_0=t_0$.
Given $s_i$, let $\xi_i\leq \eta_0$ be the greatest such that $s_i=t_{\xi_i}$, i.e., the greatest $\xi$ such that $s_i\leq_{\xi} s$.
We notice that such a greatest ordinal exists by the continuity of $(\leq_\xi)_{\xi\leq\eta}$.
(To find $\xi_i$ computably, first check if $s_i\leq_{\eta_0} s$.
If so let $\xi_i=\eta_0$, and if not search for $\xi$ with $s_i\leq_\xi s$ and $s_i\not\leq_{\xi+1} s$.)
If $\xi_i=\eta_0$, then we let $k=i$ and that finishes the definition of $s_0,...,s_k$.
Otherwise, let $s_{i+1}=t_{\xi_i+1}$.
Since we know $s_i\not\leq_{\xi_i+1}s$, we must have $s_{i+1}<s_i$.
By $(\clubsuit)$ we then have $s_{i+1}\leq_{\xi_i+1}s_i$, and hence, since $p$ is a 0-run, $\ell_{s_{i+1}}\leq_{\xi_i+1}\ell_{s_i}$.
(See diagram below.)
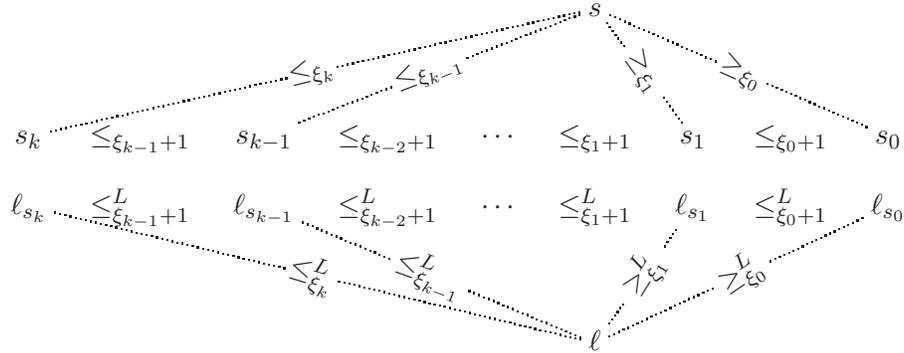
\[  \label{big diagram}
\xymatrix@C=10pt@R=7pt{
	& 			& 		& 			&		& 	 
	s \ar@{.}|{\rotatebox[origin=c]{14}{$\leq_{\xi_k}$}}[lllllddd] \ar@{.}|{\rotatebox[origin=c]{21}{$\leq_{\xi_{k-1}}$}}[lllddd]  \ar@{.}|{\rotatebox[origin=c]{-54}{$\geq_{\xi_{1}}$}}[rddd]  \ar@{.}|{\rotatebox[origin=c]{-24}{$\geq_{\xi_0}$}}[rrrddd]  \\ \\ \\ 
{{s_k}}  &  \leq_{\xi_{k-1}+1} &   {s_{k-1}} 	& \leq_{\xi_{k-2}+1} &  \cdots &  \leq_{\xi_1+1}  &  {{s_1}}  &  \leq_{\xi_0+1}    &  {{s_0}}  \\ 	
{\ell_{s_k}}  & \leq^L_{\xi_{k-1}+1} &   \ell_{s_{k-1}} 	& \leq^L_{\xi_{k-2}+1} &  \cdots &  \leq^L_{\xi_1+1}   &  {\ell_{s_1}}  & \leq^L_{\xi_0+1}    &  {\ell_{s_0}}  \\ 	\\
\\
	& 			& 		& 			&		& 	  
	{\ell} \ar@{.}|{\rotatebox[origin=c]{-17}{$\leq^L_{\xi_k}$}}[llllluuu] \ar@{.}|{\rotatebox[origin=c]{-20}{$\leq^L_{\xi_{k-1}}$}}[llluuu]  \ar@{.}|{\rotatebox[origin=c]{58}{$\geq^L_{\xi_{1}}$}}[ruuu]  \ar@{.}|{\rotatebox[origin=c]{34}{$\geq^L_{\xi_0}$}}[rrruuu] 
}\]
So we are in the hypothesis of the extendibility condition, which gives us an $\ell\in L$ with $p\concat\ell\in P_0$ and $\ell_{s_{i}}\leq_{\xi_i}\ell$ for all $i\leq k$.
As we explained above, this is exactly what we needed.
\end{proof}

%%%%%%%%%%%%%%%%%%%%%%%%%%%%%%%%
\section{The modified $\eta$-system}  \label{se: small modification}

Just by looking at the proof of our metatheorem (Theorem \ref{thm: meta eta}), one can see that we could have used the following weakening of the extendibility condition to get the same result.
The only difference is that we require $s_i\leq_{\xi_i}s$ in the hypothesis. 

\begin{definition}
We say that an $\eta$-system $(L, P_0, (\leq^L_\xi)_{\xi\leq\eta})$ satisfies the {\em weak extendibility condition} if:
whenever we have a finite 0-run $p=\la \ell_0,...,\ell_{s-1}\ra$, stages $s_k < s_{k-1} < ....< s_0< s$, and ordinals $\xi_0<\xi_1<\cdots<\xi_k\leq\eta$ such that, for all $i\leq k$, $s_i \leq_{\xi_i} s$ and $\ell_{s_{i+1}} \leq^L_{\xi_i+1} \ell_{s_i}$, there exists an $\ell\in L$ such that $p\concat\ell\in P_0$ and, for all $i\leq k$, $\ell_{s_i} \leq^L_{\xi_i} \ell$.
(See the diagram inside the proof of Theorem \ref{thm: meta eta}.)
\end{definition}

The assumption that $s_i \leq_{\xi_i} s$ for all $i\leq k$ (that is, the top half of the diagram) is not part of Ash's original framework.
Ash's framework is equivalent to the one we described in the previous section, and there is no clear way of expressing the  weak extendibility condition using Ash's original exposition as he does not deal with the $\xi$-belief relations.
The weak extendibility condition makes the metatheorem stronger (or at least so it seems), and it is necessary for the application we provide in this paper.

%The advantage of this assumption is that we know that for each $i$, $p\upto_{l_{s_i}+1}\in P_0$, and this extra-information about the $\ell_{s_i}$'s is important in our application for showing that $\ell$ exists.
%As witnessed by all the applications that appear in Ash and Knight's book, this assumption if often not necessary.

\begin{theorem}\label{thm: meta eta weak}
For every $\eta$-system $(L, P_0, (\leq^L_i)_{i\leq \eta})$ with the weak extendibility condition, there exists a computable infinite 0-run $\pi$.
\end{theorem}

The proof of this theorem is exactly that of Theorem \ref{thm: meta eta}.

%%%%%%%%%%%%%%%%%%%%%%%%%%%%%%%%%%%%%%%%%%%%%%%%%%%%%%%%%%%%%%%%%%%%%%%%%%%%%%%%%%%%%%%%%%%%%%%%%%%%%%%%%%%%%%%%%%%%%%%%%%%%%%%%%%%%%%%%%%%%%%%%%%%%%%%%%%%%%%%%%%%%%%%%%%%%%%%%%%%%%%%%%%%%%%%%%%%%%%%%%%%%%%%%%%%%%%%%%%%%%%%%%%%%%%%%%%%%%%%%%%%%%%%%%%%%%%%%
\section{An application}  \label{se: app}

In this section we prove Theorem \ref{thm: app} which requires the use of an $\eta$-system with the weak extendibility condition.
This theorem is new, and is key in the author's upcoming paper \cite{MonIntermediateVC}.
The proof is a modification of the proof of Ash and Knight's pairs-of-structures theorem \cite[Theorem 18.6]{AK90}.
To state our result we need a couple definitions. 
Fix a computable ordinal $\eta$, and suppose that the definitions of $\nabla^{\xi}_s$ and $\leq_\xi$ hold for all $\xi\leq \eta$.

\begin{definition}
Let $\fseq{\eta}$  be the set of all $\si\in 2^\eta$ with only finitely many 1's.
\end{definition}

Notice that $\fseq{\eta}$ is countable and computably presentable, as opposed to $2^\eta$ which has size continuum for infinite $\eta$. 

Let us recall the back-and-forth relations.
For more background on the back-and-forth relation see \cite[Chapter 15]{AK00}.
Given structures $\A$ and $\B$, and tuples $\abar\in\A^{<\om}$ and $\bbar\in\B^{<\om}$, recall that $(\A,\abar)$ is $\xi$-back-and-forth below $(\B,\bbar)$, written $(\A,\abar)\leq_\xi(\B,\bbar)$, if the $\Pii_\xi$-type of $\abar$ in $\A$ is included in the $\Pii_\xi$-type of $\bbar$ of $\B$, where the {\em $\Pii_\xi$-type of $\abar$ in $\A$} is the set of infinitary-$\Pi_\xi$ formulas which are true about $\abar$ in $\A$.
(We are allowing tuples of different sizes here as in \cite{AK00}, provided $|\abar|\leq |\bbar|$.
We note that $(\A,\abar)\leq_\xi(\B,\bbar) \iff (\A,\abar)\leq_\xi(\B,\bbar\upto |\abar|)$.)
Equivalently, $(\A,\abar)\leq_\xi(\B,\bbar) $ if for every tuple $\dbar\in\B^{<\om}$ and any $\g<\xi$, there exists $\cbar\in\A^{<\om}$ such that $(\A,\abar\cbar)\geq_\g(\B,\bbar\dbar)$.
In particular, we will use that, if $(\A,\abar)\leq_{\b+1}(\B,\bbar)$, then there exists a tuple $\cbar\supseteq\abar$ such that  $(\A,\cbar)\geq_{\b}(\B,\bbar)$.
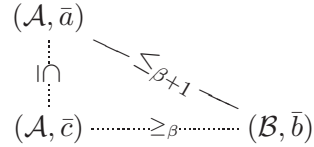
\[
\xymatrix{ 
(\A,\abar)\\
  (\A,\cbar)  \ar@{.}|{\rotatebox[origin=c]{-90}{$\subseteq$}}[u]    \ar@{.}|{\geq_{\b}}[rr]     & &      (\B,\bbar)  \ar@{-}|{\rotatebox[origin=c]{-30}{$\leq_{\b+1}$}}[llu] 
}\]

\begin{definition}
An {\em $\eta$-tree of structures} is a sequence of structures $\{\A_\si: \si\in \fseq{\eta}\}$ such that, for every $\si,\tau\in \fseq{\eta}$ and $\xi<\eta$, we have that 
\[
\si\upto \xi=\tau\upto \xi \quad\and\quad \si(\xi)\leq\tau(\xi)
\quad \implies \quad
\A_{\si}\geq_{\xi+1}\A_{\tau}.
\]
\end{definition}
Notice that, in particular $\si\upto \xi=\tau\upto \xi \implies \A_{\si}\equiv_{\xi}\A_{\tau}$.

We now recall the notion of $\a$-friendliness (see \cite[Section 15.2]{AK00}).
An $\eta$-tree of structures is {\em $\a$-friendly} if given two tuples in two structures $\abar\in\A_\si^{<\om}$ and $\bbar\in\A_\tau^{<\om}$ and given $\xi<\a$, we can effectively decide if $(\A_\si,\abar)\leq_\xi(\A_\tau,\bbar)$ in a c.e.\ way, or, in other words, if the set of quintuples $\{(\si,\abar,\tau,\bbar,\xi): \si,\tau\in \fseq{\eta},\abar\in\A_\si^{<\om}, \bbar\in\A_\tau^{<\om}, \xi<\a, (\A_\si,\abar)\leq_\xi(\A_\tau,\bbar)\}$ is c.e.
This is quite a strong condition, but when the $\A_\si$'s are naturally defined structures which we fully understand, it is often the case that we can show the $\eta$-tree is $\a$-friendly.

The following is our generalization of a Ash--Knight's theorem from pairs to trees of structures.

\begin{theorem}\label{thm: app}
Let $\{\A_\si:\si\in\fseq{\eta}\}$ be a computable $\eta$-friendly $\eta$-tree of structures.
Let $\{\si_n:n\in\om\}\subseteq \fseq{\eta}$ be such that deciding if $\si_n(\xi)=1$ is $\Sigma^0_{\xi+1}$ uniformly in $n$ and $\xi$.

Then, there exists a computable sequence $\C_n$ such that for all $n$, $\C_n\isom\A_{\si_n}$.
\end{theorem}

In other words, if we are given a $\si\in\fseq{\eta}$, and it is given to as in a $\Sigma^0_{\xi+1}$ way, we can uniformly build a computable copy of $\A_\si$, despite having to guess at the bits of $\si$.

\begin{proof}
To build $\C_n$ we will define an $\eta$-system $(L, P_0, (\leq^L_i)_{i\leq \eta})$ uniformly in $n$.
We will have to change the definition of ``apparent $\xi$-true stage'' slightly.
That is, we will define a new sequence of pre-orderings $(\Aleq_\xi)_{\xi\leq\eta}$ satisfying (\ref{pa: 0})-(\ref{pa:6}) and $(\clubsuit)$ so that we can still apply the metatheorem.

Let $W$ be a c.e.\ operator such that $\si_n(\xi)=1\iff n\in W^{\nabla^{\xi+1}}$.
Let us fix an $n$ and concentrate on building $\C_n$.

For each $s$, let $\tau_s\in\fseq{\eta}$, the stage-$s$ approximation to $\si_n$, be defined as follows:
\[
\tau_s(\xi)=1\iff n\in W^{\nabla^{\xi+1}_s}.
\]
Since $\nabla^{\xi+1}_s\neq\la\ra$ only for finitely many $\xi$'s, we have that $\tau_s(\xi)=1$ also only for finitely many $\xi$'s.
We note that if $\nabla^\eta_t$ is correct, and $t$ is large enough so that $n\in W^{\nabla^{\xi+1}_t}$ for all $\xi$ with $\si_n(\xi)=1$, then since all the $\nabla^{\g}_t$ are correct, $\tau_s=\sigma_n$.

For $s,t\in\om$ let us define
\[
s \Aleq_\xi t   \quad\iff\quad s\leq_\xi t \ \ \and \ \ \tau_s\upto\xi=\tau_t\upto\xi.
\]

First, let us note that $s \Aleq_\xi t$ not only implies that $\tau_s\upto\xi=\tau_t\upto\xi$, but also that $\tau_s(\xi)\leq\tau_t(\xi)$:
This is because if $s \Aleq_\xi t$, then $s \leq_\xi t$, and hence $\nabla^{\xi+1}_s\subseteq\nabla^{\xi+1}_t$, and thus $n\in W^{\nabla^{\xi+1}_s} \implies n\in W^{\nabla^{\xi+1}_t}$.

It is easy to see that the $(\Aleq_\xi)_{\xi\leq\eta}$ satisfy (\ref{pa: 0})-(\ref{pa:6}).
To show $(\clubsuit)$ suppose that $r\Aleq_\xi s\Aleq_\xi t$ and $r\Aleq_{\xi+1}t$--we want to show that $r\Aleq_{\xi+1}s$.
Since $\leq_\xi$ satisfies $(\clubsuit)$, we have that $r\leq_{\xi+1}s$.
Since $r\Aleq_\xi s$ we have that $\tau_r\upto\xi=\tau_s\upto\xi$.
Since $r\Nleq_\xi s\Nleq_\xi t$,  $\tau_r(\xi)\leq \tau_s(\xi)\leq\tau_t(\xi)$ (as we mentioned above) and since $r\Aleq_{\xi+1}t$, $\tau_r(\xi)= \tau_t(\xi)$.
It follows that $\tau_r(\xi)=\tau_s(\xi)$, getting that $\tau_s\upto\xi+1=\tau_t\upto\xi+1$ as needed.

We will use $(\Aleq_\xi)_{\xi\leq\eta}$ as our ``$\xi$-belief'' relations.
Let us now define the $\eta$-system.
At each stage $s$ we will pick a tuple $\abar_s\in A_{\tau_s}$ and build a piece of our structure $\C_n$ using the atomic diagram of $\abar_s$.

We let 
\[
L=\{ (\tau,\abar): \tau\in \fseq\eta, \abar\in A_\tau\}.
\]

We only play $(\tau,\abar)$ when the current guess at $\si_n$ is $\tau$, and we want $\abar$ to cover more and more of $\A_\tau$ every time.
We define the action tree accordingly:
 $p\concat(\tau,\abar)\in P_0$ if and only if $p\in P_0$, $\tau=\tau_{s}$, where $s=|p|$, and $\abar$ contains the first $s$ elements of $\A_\tau$.

Let $c_0,c_1,....$ be a fresh set of constants.
We will define $\C_n$ to have domain $\{c_0,c_1,....\}$, (quotiented by the equality relation).
To define the structure we need to define the atomic formulas among these constants.
The enumeration function will enumerate the atomic type of the tuple we are currently considering.
That is,  given a literal $\psi(x_{i_1},...,x_{i_k})$, we let $\psi(c_{i_1},...,c_{i_k})\in E(\tau,\abar)$ if the G\"odel code of $\psi$ is less than $|\abar|$, and if $\A_\tau\models\psi(a_{i_1},...,a_{i_k})$ where $\abar=\la a_0,...,a_{|\abar|-1}\ra$.
We will let $(\tau,\abar)\leq^L_0 (\rho,\bbar)$ if $(\A_\tau,\abar)\leq_0(\A_\rho,\bbar)$, which implies $E(\tau,\abar)\subseteq E(\rho,\bbar)$.
So, if we have a computable infinite 0-run $\pi$, $E(\pi)$ is the atomic diagram of a computable structure.

We define the top restraint relation $\leq^L_\eta$ as follows
\[
(\tau_0,\abar_0)\leq^L_\eta(\tau_1,\abar_1) 
	\quad\iff\quad
\tau_0=\tau_1 \ \ \and\ \ \abar_0\subseteq\abar_1.
\]
Thus, if we have an infinite 0-run $\pi$, and $t_0\Aleq_\eta t_1\Aleq_\eta \cdots$ is an infinite sequence given by (\ref{pa:6}), then $\tau_{t_i}=\si_n$ and $\abar_{t_i}\subseteq \abar_{t_{i+1}}$ for all $i$.
Therefore, putting together these tuples we get a function $\vec{a}=\bigcup_i{\abar_{t_i}}\colon\om\to\A_{\si_n}$.
From the definition of $E(\pi)$ and $\C_n$, we notice that the map $c_k\mapsto \vec{a}(k)$ is then an isomorphism between $\C_n$ and $\A_{\si_n}$.
So, if we find a computable infinite 0-run, we have a computable $\C_n\isom\A_{\si_n}$ as wanted.

To finish the definition of the $\eta$-system we need to define the intermediate pre-orderings.
For $0<\xi<\eta$ we let $(\tau_0,\abar_0)\leq^L_\xi(\tau_1,\abar_1)$ if the appropriate back-and-forth relation holds, that is, $(\A_{\tau_0},\abar_0)\leq_\xi(\A_{\tau_1},\abar_1)$.

To show that an infinite computable run exists and finish the proof, we need to prove the weak extendability condition.
Suppose we have a finite 0-run $\la \ell_0,...,\ell_{s-1}\ra$, a sequence of stages $s_k<s_{k-1}<...<s_0< s$, and a sequence of ordinals $\xi_0<\xi_1<...<\xi_k\leq\eta$ such that, for all $i$,
\[
s_i\Aleq_{\xi_i} s
\quad \and\quad
\ell_{s_{i+1}} \leq^L_{\xi_{i}+1} \ell_{s_{i}}.
\]

For each $t<s$, let $\la\tau_s,\abar_s\ra=\ell_s$.
The first step is to fix the mistakes made by all the stages between $s_k$ and $s$.
By recursion on $j\leq k$ we build a tuple $\bbar_j\in\A_{\tau_{s_j}}$ such that $\abar_{s_j}\subseteq\bbar_j$ and such that for all $i\leq j$, $\ell_{s_i}\leq^L_{\xi_i}(\tau_{s_j},\bbar_j)$.
(See diagram below.)
Let $\bbar_0=\abar_{s_0}$.
Having built $\bbar_j$, since $\abar_{s_j}\subseteq\bbar_j$, we know that $(\A_{\tau_{s_{j+1}}},\abar_{s_{j+1}})\leq_{\xi_j+1}(\A_{\tau_{s_{j}}},\bbar_j)$.
But then, by the basic back-and-forth property, there exists $\bbar_{j+1}\supseteq\abar_{s_{j+1}}$ such that  $(\A_{\tau_{s_{j+1}}},\bbar_{j+1})\geq_{\xi_j}(\A_{\tau_{s_{j}}},\bbar_j)$, and hence for each $i\leq j$, $(\A_{\tau_{s_{j+1}}},\bbar_{j+1})\geq_{\xi_i}(\A_{\tau_{s_{i}}},\abar_{s_i})$.

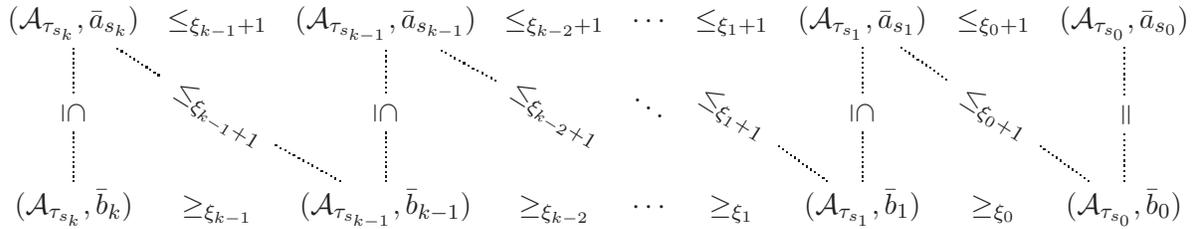
\[
\xymatrix@C=4pt@R=7pt{
(\A_{\tau_{s_k}},\abar_{s_k})    \ar@{.}[d]\ar@{.}[rd]   		&\leq_{\xi_{k-1}+1}&   (\A_{\tau_{s_{k-1}}},\abar_{s_{k-1}})  \ar@{.}[d]\ar@{.}[rd]   
		&\leq_{\xi_{k-2}+1} &  \cdots &  \leq_{\xi_{1}+1}   &  
			(\A_{\tau_{s_1}},\abar_{s_1})  \ar@{.}[d]\ar@{.}[rd]    & \leq_{\xi_0+1}    &  (\A_{\tau_{s_0}},\abar_{s_0})   \ar@{.}[d]  \\ 	
\rotatebox[origin=c]{-90}{$\subseteq$} & \rotatebox[origin=c]{-32}{$\leq_{\xi_{k-1}+1}$} & \rotatebox[origin=c]{-90}{$\subseteq$} & \rotatebox[origin=c]{-32}{$\leq_{\xi_{k-2}+1}$} &  \ddots  & \rotatebox[origin=c]{-32}{$\leq_{\xi_{1}+1}$} & \rotatebox[origin=c]{-90}{$\subseteq$} & \rotatebox[origin=c]{-32}{$\leq_{\xi_{0}+1}$} & \rotatebox[origin=c]{-90}{$=$} \\
(\A_{\tau_{s_k}},\bbar_{k})   \ar@{.}[u]    & \geq_{\xi_{k-1}}    &  (\A_{\tau_{s_{k-1}}},\bbar_{{k-1}})  \ar@{.}[u]\ar@{.}[lu] 	&  \geq_{\xi_{k-2}}  		&\cdots      &   \geq_{\xi_{1}}   & (\A_{\tau_{s_{1}}},\bbar_{1}) \ar@{.}[u]\ar@{.}[lu]&   \geq_{\xi_0}&     (\A_{\tau_{s_0}},\bbar_{0}) \ar@{.}[u]   \ar@{.}[lu] 
}\]

We note that if it was the case that $\xi_{k-1}+1=\xi_k=\eta$, we would have $\ell_{s_k}\Aleq_\eta\ell_{s_{k-1}}$ and hence that $\abar_{s_k}\subseteq\abar_{s_{k-1}}$, so we could just let $\bbar_k=\bbar_{k-1}$.

The second step is to define $\cbar\in \A_{\tau_{s}}$ such that  $(\A_{\tau_{s_k}},\bbar_k)\leq_{\xi_k} (\A_{\tau_{s}},\cbar)$.
(This is the step that requires the weak version of the extendibility condition.)
Since $s_k\Aleq_{\xi_k} s$ we have that $\tau_{s_k}\upto \xi_k=\tau_{s}\upto \xi_k$ and  $\tau_{s_k}(\xi_k)\leq \tau_{s}(\xi_k)$.
Thus, from the definition of $\eta$-tree we have that $\A_{\tau_{s_k}}\geq_{\xi_k+1} \A_{\tau_{s}}$.
By the back-and-forth property, there exist $\cbar\in \A_{\tau_{s}}$ such that $(\A_{\tau_{s_k}},\bbar_k)\leq_{\xi_k} (\A_{\tau_{s}},\cbar)$ exactly as wanted.

Last, we  define $\abar\supseteq \cbar$ so that $p\concat (\tau_{s},\abar)\in P_0$:
All we need to do is extend $\cbar$ making sure the first $s$ elements of $\A_{\tau_{s}}$ are included in it.
\end{proof}

%%%%%%%%%%%%%%%%%%%%%%%%%%%%%%%%%%%%%%%%%%
\section{Approximation of the $\xi$-jumps} \label{se: true stages}

The objective of this section and the next is to define, for each ordinal $\xi$,  the notion of ``$\xi$-true stage'' and the one of ``apparent $\xi$-true stage'' as described in Section \ref{see: apparent true}.
To get these notions to work smoothly we need to make a few technical definitions and modify a few known notions so that they all fit together.
In this section we concentrate on defining the approximations to the $\xi$-jumps.

Most of the notions below were defined by Marcone and Montalb\'an in \cite{MMVeblen}.
In that paper a whole lot of properties are proved in much more detail than here.
In that paper we were interested in comparing the behavior of the iterates of the Turing jump and certain ordinal-notation functions that come from proof theory.
The reader may just look at Sections 4, 5.2 and 6.1 in \cite{MMVeblen} and skip the proof-theoretic parts.
We will repeat all the relevant definitions and properties below.

As usual, let $\varphi_0,\varphi_1,...$ be an effective enumeration of the Turing functionals.
%Changing what functional is represented by $\varphi_0$ will be useful in applications, but it is, of course, irrelevant in this section and the next.

%%%%%%%%%%%%%%%%%%%%%%%
\subsection{Approximation of the first jump}

The idea of {\em  true stage in the enumeration of $0'$} is widely used in Computability Theory, and there are a few slightly different definitions in the literature.
We use our own.
We want the $i$-th true stage to be a stage in the enumeration of $0'$ where the first $i$ bits have been guessed correctly.
For technical reasons, we will also require the $i$-th true stage to be larger than the $(i-1)$-th true stage.
Also, we want to be able to iterate and relativize this definition.

\begin{definition}\label{def:1true}
Given  $Z \in \Bai$, we define the sequence of {\em $Z$-true stages} as follows:
\[
t_{i} = \begin{cases}   t_{i-1}+1,		        									& \mbox{if }\varphi_{i}^{Z}(i)\diverges,  \\
				\max\{t_{i-1}+1, \mu t(\varphi_{i}^{Z\upto t}(i)\converges)\},   	& \mbox{if }\varphi_{i}^{Z}(i)\converges,  \end{cases}
\]
starting with $t_{-1}=1$ (and hence with $t_0>1$). 
\end{definition}

So, $t_i$ is a stage where $Z$ can correctly guess
$Z'\upto i+1$ because 
\[
\forall m \leq i (m\in Z'\iff \varphi_m^{Z\upto t_i}(m) \converges).
\] 

(We use the standard convention that if $\si$ is a finite string, then $\varphi_e^\si(i)$ runs for at most $|\si|$ stages.)
With this in mind, we give a new definition of the Turing jump operator.

\begin{definition}\cite[4.1]{MMVeblen}  \label{def:J}
We define the \emph{Jump operator} to be the function $\J\colon \Bai \to \Bai$ such that for every $Z \in \Bai$ and $n \in \om$,
\[
\J(Z)(n) = \ulcorner Z\upto t_n\urcorner,
\]
by which we mean the natural number coding the string $Z\upto t_n$.
(We will abuse notation a write $Z\upto t_n$ for both the string $Z\upto t_n$ and the number $\ulcorner Z\upto t_n\urcorner$ depending on the context.)
Equivalently
\[
\J(Z)=\la Z\upto t_0, Z\upto t_1, Z\upto t_2, Z\upto t_3, \dots\ra.
\]
\end{definition}

Clearly $\J(Z) \teq Z\oplus \la t_0,t_1,t_2,...\ra \teq Z'$.
So we will use $\J(Z)$ as our standard jump operation because it plays nicely with the approximations for the jump that we define below.

\begin{definition}  \cite[4.2]{MMVeblen}  \label{def: J si}
The \emph{jump approximation function} is the map $J \colon \Seq \to \Seq$ defined as
follows. For $\si \in \Seq$, define 
\[
t_{i} = \begin{cases}  t_{i-1}+1,		        									& \mbox{if }\varphi_{i}^{\si}(i)\diverges,  \\
				\max\{ t_{i-1}+1, \mu t(\varphi_{i}^{\si\upto t}(i)\converges)\},   	& \mbox{if }\varphi_{i}^{\si}(i)\converges,  \end{cases}
\]
starting with $t_{-1}=1$. 
Let 
\[
J(\si)=\la\si\upto t_0, \si\upto t_1, \dots,\si\upto t_{k}\ra
\]
where $k$ is the least such that $t_{k+1}>|\si|$. So, either $|\si|\leq 1$ and $J(\si)=\la\ra$, or $k$ is the first with $t_k=|\si|$.
\end{definition}

The reason we want to start with $t_0\geq 2$ is that we then have $|J(\si)|<|\si|$.
We remark that the jump approximation function is computable.
Another important remark is that $J(Z\upto t)\subset\J(Z)$ if and only if $t$ is a $Z$-true stage.

A nice combinatorial property $J$ has is that its inverse function is very easy to compute, as $\si$ is coded as the last entry of $J(\si)$.
Given $\tau\in\Seq$, we let $J^{-1}(\tau)\in\Seq$ be the string coded by the last entry of $\tau$.
(In \cite{MMVeblen} we used the letter $K$ to denote the map $J^{-1}$.)
So, for $|\si|>1$,  $J^{-1}(J(\si))=\si$.

The following simple lemma is the key to get the $(\clubsuit)$ property later.

\begin{lemma} \label{le: s t u}
If $\si\subseteq \tau \subseteq \pi$ and $J(\si) \subseteq J(\pi)$, then $J(\si) \subseteq J(\tau)$.
\end{lemma}
\begin{proof}
The reason is that, if the values of $t_0,...,t_k$ as in the definition of $J(\si)$ remind unchanged in the definition of $J(\pi)$ (except maybe getting a larger $k$), it is because for $i\leq k$, if $\varphi_i^\si(i)\diverges$, then $\varphi_i^\pi(i)\diverges$ too. 
But then, $\varphi_i^\tau(i)\diverges$ also, and the values of $t_0,...,t_k$ remain unchanged in the definition of $J(\tau)$ as well.
\end{proof}

\begin{lemma} \label{le: J properties}
If  $\la\ra\neq J(\si)\subseteq J(\tau)$, then $\si \subseteq \tau$.
\end{lemma}
\begin{proof}
This is proved in \cite[Lemma 4.4 (P2)(P6)]{MMVeblen}.
The reason is that the last entry of $J(\si)$ codes the string $\si$ itself, and any entry of $J(\tau)$ codes an initial segment of $\tau$.
\end{proof}

%%%%%%%%%%%%%%%%%%
\subsection{The $\om$-jump}

For finite $n$, the $n$-th jump and its approximation are just defined by iterating $\J$ and $J$.
That is, we let $\J^{n}=\J\circ\J\circ\cdots\circ\J\colon\Bai\to\Bai$ ($n$ times) and $J^{n}=J\circ J\circ\cdots\circ J$ ($n$ times).
We still get that, for each $Z\in\Bai$, there is a sequence of $n$-$Z$-true stages $\{t^n_i:i\in\om\}$ such that  $\J^n(Z)=\bigcup_{i\in\om}J^n(Z\upto t^n_i)$, and everything works quite smoothly.
We can also define $J^{-n}=J^{-1}\circ\cdots\circ J^{-1}$ ($n$ times), and we still get that $J^{-n}\circ J^n$ is the identity on stings $\si$ with $J^n(\si)\neq\la\ra$.
It is not hard to see that $|J^n(\si)| \leq |\si|\dmi n$, where $s\dmi t=\max(0,s-t)$.

The reader only interested in $n$-systems for finite $n$ may skip the rest of this section.

The $\om$-th jump is usually defined by putting all the $n$ jumps together, each in a different column.
Instead, we will take just one entry from each.
This is not going to change the Turing degree, but it is going to play better with our finite approximations. 

\begin{definition}
We define the \emph{$\om$-Jump operator} to be the function
$\J^\om\colon \Bai \to \Bai$ such that for every $Z \in \Bai$ and $n\in\om$
\[
\J^\om(Z)(n)=\J^{n+1}(Z)(0),
\]
or, in other words, 
\[
\J^\om(Z) = \la \J^1(Z)(0),\  \J^2(Z)(0),\  \J^3(Z)(0),\dots\ra,
\]
\end{definition}

Not unexpected is that $\J^\om(Z) \teq \bigoplus_n\J^n(Z) \teq Z^{(\om)}$ for every $Z\in \Bai$.
The proof is exactly as the one given in \cite[Lemma 5.7]{MMVeblen}.
The idea for the proof is that for $n<m$, $J^{-(m-n)}(\la\J^m(Z)(0)\ra)$ is an initial segment of $\J^n(Z)$ of length at least $m-n$.
Thus $\J^n(Z)$ can be uniformly defined from $\J^{\om}(Z)$ using that $\J^{n}(Z)=\bigcup_{m\in\om}J^{-(m-n)}(\la\J^m(Z)(0)\ra)$.

\begin{definition}\label{def:Jom si}
The \emph{$\om$-jump approximation function} is the map $J^\om \colon \Seq \to \Seq$
defined as follows. Given $\si\in \Seq$, let
\[
J^\om(\si) = \la J^1(\si)(0),\  J^2(\si)(0),\dots,  J^{n}(\si)(0)\ra,
\]
where $n$ is the least such that $J^{n+1}(\si)=\emptyset$. 
(There is always such an $n$, because $|J^i(\si)|\leq|\si|\dmi i$.) 
\end{definition}

%This definition deviates slightly from that in \cite[Definition 5.4]{MMVeblen} in that we start with $J^2(Z)(0)$, rather than with $\J(Z)(0)$.
%The reason we do this is that we want to have $|J^\om(\si)|\leq|\si|\dmi 2$.
%To see this note that if $|J^\om(\si)|=k$, then $|J^{k+1}(\si)|\geq 1$, and hence $|\si|\geq k+2$.

The inverse of $J^\om$ is still easy to compute.
We let $J^{-\om}(\tau) = J^{-|\tau|}(\la last(\tau)\ra)$.
It is not hard to see that $J^{-\om}$ is the inverse of $J^\om$ on strings $\si$ with $J^\om(\si)\neq\la\ra$.

\begin{lemma} \label{le: omega inclusion}
For $\si,\tau\in\Seq$, $J^\om(\si)\subseteq J^\om(\tau)$ if and only if $(\forall n<\om)\ J^n(\si)\subseteq J^n(\tau)$, provided $\la\ra\neq J^\om(\si)$.
\end{lemma}
\begin{proof}
The right-to-left direction is clear, as $J^n(\si)\subseteq J^n(\tau)$ implies $J^n(\si)(0)= J^n(\tau)(0)$ whenever defined.
For the left-to-right direction we have that if $|J^\om(\si)|=k$, then $|J^{k}(\si)|=1$ and hence $J^{k}(\si)=\la J^\om(\si)(k-1)\ra=\la J^\om(\tau)(k-1)\ra= \la J^{k}(\tau)(0)\ra\subseteq J^{k}(\tau)$.
For each $j\leq k$, applying Lemma \ref{le: J properties} $k-j$ times we get that $J^j(\si)\subseteq J^j(\tau)$.
For $j>k$ we have $J^j(\si)=\la\ra$.
\end{proof}

%%%%%%%%%%%
\subsection{The transfinite jumps}

For general $\g$ the definitions of the $\om^\g$-th jump operator and the $\om^\g$-jump approximation function are slightly more involved, but only because of the notation, while the ideas are the same as in the $\om$-jump.

Assume again that $\eta$ is a computable ordinal.
From now on we only work with ordinals in $\eta+1$.
Since the jump operators are going to be defined by transfinite recursion, we need an effective way to reach each limit ordinal from below.
Let us review the notion of characteristic sequence of an ordinal $\a$.
It is just a non-decreasing sequence of ordinals, that we denote by $\{\a[0],\a[1],...\}$, which satisfies  $\lim_n\a[n]+1=\a$, and that we fix in advance.
Each ordinal $\a$ is going to come together with such a sequence, (like when we have an ordinal notation for the ordinal $\a$). 
It is not hard to see that such a sequence can be found computably, uniformly in $\a$.
For instance, if  $\a$ is a successor ordinal, we might have $\a[n]=\a-1$ for all $n$.
%If $\a$ is a limit ordinal we define it characteristic sequence so that $\a[0]<\a[1]<\cdots$ is a sequence with limit $\a$.
If $\a=0$, then $\a[n]$ is undefined.
%We ask that if $\a$ and $\b$ are different limit ordinals, then $\a[0]\neq\b[0]$.
%Finding such sequence computable is straightforward.
%For the injectivity property note that we can always chose $\a[0]$ such that, if $\b$ is a limit ordinal with $\ulcorner\b\urcorner <\ulcorner \a\urcorner$, then $\a[0]\neq\b[0]$.
%and $\b[0]<\a<\b$, we chose $\a[0]$ so that $\b[0]<\a[0]<\a$.
%So, in particular, we get that each $\g$ is use as $\a[0]$ for at most one limit $\a$.
Observe that in both the limit and the successor cases we have that $\om^\a=\sum_{i\in\om}\om^{\a[i]}$.

%First, if $\g$ is not of the form $\om^\a$, then we write $\g$ as $\om^{\a_0}+\om^{\a_1}+\cdots+\om^{\a_k}$ with $\a_0\geq\cdots\geq\a_k$, and define
%\[
%\J^\g=\J^{\om^{\a_k}}\circ\cdots\circ\J^{\om^{\a_0}}\colon\Bai\to\Bai    \quad \and \quad   J^\g=J^{\om^{\a_k}}\circ\cdots\circ J^{\om^{\a_0}}\colon\Seq\to\Seq.
%\]
%The inverse function is easily definable as $J^{-n}\g=J^{-\om^{\a_0}}\circ\cdots\circ J^{-\om^{\a_k}}$, and is clear that whenever $J^\g(\si)\neq\la\ra$, $J^{-n}\g(J^\g(\si))=\si$.

To simplify the notation in the definition of the $\om^\a$-Jump function, let $\J^{\om^\a}_n\colon \Bai \to \Bai$ and
$J^{\om^\a}_{ n}\colon \Seq \to \Seq$ be defined by
\begin{align*}
\J^{\om^\a}_{ n} &= \J^{\om^{\a{[n-1]}}}\circ \J^{\om^{\a[n-2]}} \circ \cdots \circ \J^{\om^{a[0]}}, \\
\Jom\a_{ n} &= J^{\om^{\a[n-1]}}\circ J^{\om^{\a[n-2]}} \circ \cdots \circ \Jom{\a[0]}.
\end{align*}

Note that $\J^{\om^\a}_n(Z)\teq Z^{(\b_n)}$ where $\b_n=\om^{\a[0]}+\cdots+\om^{\a[n-1]}$ and that $\lim_n\b_n=\om^\a$.
It thus follows that $Z^{(\om^\a)}\teq \bigoplus_n Z^{(\b_n)}\teq \bigoplus_n \J^{\om^\a}_n(Z)$.
The definition of $\J^{\om^\a}(Z)$ will, instead, use just one bit--the first bit--from each $\J^{\om^\a}_n(Z)$.

%If $\g=\om^\a$, for each $n$, let $\b_n=\sum_{i=0}^{n}\om^{\a[i]}$, so that $\om^\a=\lim_n\b_n$.
\begin{definition}(\cite[Definitions 6.1 and 6.5 ]{MMVeblen})
Define $\J^{\om^\a}$ by transfinite recursion as follows.
For $Z\in\Bai$ and $n\in\om$ let
\[
\J^{\om^\a}(Z)(n)=\J^{\om^\a}_{ n+1}(Z)(0).
\]

The \emph{$\om^\a$-jump approximation function} is the map $J^{\om^\a} \colon \Seq \to
\Seq$ defined by
\[
\Jom\a(\si) = \la \Jom\a_1(\si)(0), \Jom\a_2(\si)(0), \dots,\Jom\a_{n-1}(\si)(0)\ra,
\]
where $n$ is least such that $\Jom\a_{ n}(\si)=\emptyset$.
By computable transfinite recursion one can show that $J^{\om^\a}$ is computable.
Again, it is not hard to see that $\si$ can be recovered from the last entry of $J^{\om^\a}(\si)$:
A function $J^{-\om^\a}$ which is the inverse of $J^{\om^\a}$ on strings of length $\geq 2$ is defined in \cite[Definition 6.1]{MMVeblen}.
\end{definition}
%(\tau)=J^{-\b_{|\tau|-1}}(\la last(\tau)\ra)$, and notice that whenever $J^\g(\si)\neq\la\ra$, $J^{-n}\g(J^\g(\si))=\si$.

%This definition again deviate from \cite[Definitions 6.1 and 6.5 ]{MMVeblen} in that we start with $J^{\b_1}(\si)(0)$ rather than $J^{\b_0}(\si)(0)$.

The main properties we will use about these approximations are described in the following three lemmas.

\begin{lemma} \label{le: omega alpha inclusion}
For $\si,\tau\in\Seq$, $J^{\om^\a}(\si)\subseteq J^{\om^\a}(\tau)$ if and only if $(\forall n<\om)\ J^{\om^\a}_{ n}(\si)\subseteq J^{\om^\a}_{ n}(\tau)$, provided $\la\ra\neq J^{\om^\a}(\si)$.
\end{lemma}
\begin{proof}
The proof is essentially the same as that of Lemma \ref{le: omega inclusion}.
The right-to-left direction is straightforward, and the left-to-right direction is proved in \cite[Lemma 6.2 (P$^{\om^\a}$7)]{MMVeblen}.
%
%The right-to-left direction is clear, as $J^n(\si)\subseteq J^n(\tau)$ implies $J^n(\si)(0)= J^n(\tau)(0)$ whenever defined.
%For the left-to-right direction we have that if $|J^\om(\si)|=k$, then $|J^{k+1}(\si)|=1$ and hence $J^{k+1}(\si)=\la J^\om(\si)(k-1)\ra=\la J^\om(\tau)(k-1)\ra= \la J^{k+1}(\tau)(0)\ra\subseteq J^{k+1}(\tau)$.
%For each $j\leq k+1$, applying Lemma \ref{le: J properties} $(k-j+1)$ times we get that $J^j(\si)\subseteq J^j(\tau)$.
%For $j>k+1$ we have $J^j(\si)=\la\ra$.
\end{proof}

As a special case, we get 
\begin{corollary}\label{co: omega alpha inclusion}
If $\la\ra\neq J^{\om^\a}(\si)\subseteq J^{\om^\a}(\tau)$, then $\si\subseteq\tau$.
\end{corollary}

\begin{lemma} \label{le: inf true stages}
For every $Z$ there is a sequence $\{t_i:i\in\om\}$ such that $\bigcup_iJ^{\om^\a}(Z\upto t_i)=\J^{\om^\a}(Z)$.
\end{lemma}
\begin{proof}
This is proved in \cite[Lemma 6.9]{MMVeblen}.
Essentially, for each $i$, first prove that $J^{-\om^{\a}}(\J^{\om^\a}(Z)\upto i)\subseteq Z$, and then let $t_i=|J^{-\om^{\a}}(\J^{\om^\a}(Z)\upto i)|$. 
\end{proof}

%%%%%%%%%%%%%%%%
\subsection{Ordinal normal forms}

The material in the rest of this section is different from what is done in \cite{MMVeblen}.
In the next subsection we will define $\nabla^\xi_s$ approximating a $\Delta^0_\xi$-complete function, and we want this approximation to satisfy properties (\ref{pa:N1})-(\ref{pa:N3}).
The problem with the jump approximations we have defined so far is that for $\a$ different than $\b$, the behaviors of $J^{\om^\a}$ and $J^{\om^\b}$ are unrelated, as it depends heavily on the choice of the characteristic functions for $\a$ and $\b$ which might be completely unrelated.
This problem becomes worse when we want to define the jump along ordinals which are not of the form $\om^\a$.
The first idea would be to use standard Cantor normal forms, but this would not give us the properties we want.
In this subsection we describe a way of writing each ordinal as a sum of ordinals of the form $\om^\a$ which fits together with the characteristic functions we will use to define the jumps.

The material in this section is only relevant if the reader is interested in transfinite iterates of the jump.

For each ordinal $\eta$, let 
\[
T_\eta=\{\la n_0,...,n_k\ra\in \om^{<\om}: \eta[n_0][n_1]\cdots[n_k] \mbox{ exists}\},
\]
where by $\eta[n_0][n_1]\cdots[n_k]$ we mean the successive iteration of the characteristic functions, that is, $\eta[n_0][n_1]\cdots[n_k] = (...(\eta[n_0])[n_1]\cdots)[n_k]$.
Thus, for $\eta\neq 0$, $T_\eta$ contains a branch isomorphic to $T_{\eta[n]}$ for each $n\in\om$. 
Note that $\eta[n_0][n_1]\cdots[n_k]$ exists if and only if for all $i<k$, $\eta[n_0][n_1]\cdots[n_i]>0$.
Since $\a[n]<\a$ for all $\a$ and $n$, this tree is well-founded.
For each $\eta$, we will  define a map $g_\eta\colon T_\eta\to\om^\eta+1$ by effective transfinite recursion: Let $g_\eta(\la\ra)=\om^\eta$, and for $\si=n\concat\si^-$ let $g_\eta(\si)=\om^{\eta[0]}+\cdots\om^{\eta[n-1]}+g_{\eta[n]}(\si^-)$.
We will see below that this map is actually an isomorphism  from $(T_\eta,\leq_{KB})$ to $\om^\eta+1\sminus\{0\}$, were $\leq_{KB}$ is the Kleene-Brower ordering on $\om^{<\om}$.
(Recall that $\si\leq_{KB}\tau\iff \tau\subseteq\si \vee \left((\exists i) \si\upto i=\tau\upto i\and \si(i)<\tau(i)\right)$.)

From now on we will write $\eta\la n_0,...,n_k\ra$ for $g_\eta(\la n_0,...,n_k\ra)$.
Let us observe that this map can be defined directly as follows:
Given $\la n_0,...,n_k\ra\in T\sminus\{\la\ra\}$, we have that
\begin{eqnarray}
\eta\la n_0,...,n_k\ra & = & \sum_{i=0}^{n_0-1}\om^{\eta[i]} + (\eta[n_0])\la n_1,...,n_k\ra \\
				& = &\sum_{i=0}^{n_0-1}\om^{\eta[i]} + \cdots +  \sum_{i=0}^{n_k-1}\om^{\eta[n_0]...[n_{k-1}][i]} + \om^{\eta[n_0]...[n_{k}]} \in \om^\eta,
\end{eqnarray}
and let $\eta\la\ra=\om^\eta$.
The objective of this definition is to have a kind of normal form for writing each ordinal below $\om^\eta$ as a sum of ordinals of the form $\om^\a$.
This normal form, as opposed to the Cantor normal form, is now synchronized with the characteristic functions at the limit ordinals.

\begin{lemma} \label{le: order iso} \label{le: onto}
The map $\la n_0,...,n_k\ra\mapsto \eta\la n_0,...,n_k\ra$ is an isomorphism between $(T_\eta,<_{KB})$ and $\om^\eta+1\sminus\{0\}$. 
\end{lemma}
\begin{proof}
The proof is by induction on $\eta$.
For $\eta=0$, $T_\eta=\{\la\ra\}$,  $\om^\eta+1\sminus\{0\}=\{1\}$, and $0\la\ra=1$.
For $\eta>0$ we have that for each $n$, the map $\la n_1,...,n_k\ra\mapsto \eta[n]\la n_1,...,n_k\ra$ is an isomorphism between $(T_{\eta[n]},<_{KB})$ and $\om^{\eta[n]}+1\sminus\{\la\ra\}$. 
Notice that $T_{\eta[n]}$ equals the branch of $T_\eta$ extending $\la n\ra$, and that $\sum_n\om^{\eta[n]}=\om^\eta$.
The map $\la n_0,...,n_k\ra\mapsto \eta\la n_0,...,n_k\ra$ does nothing more than pasting all these maps together.
\end{proof}

\subsection{Special jump operators} \label{see: xi-true stages}

We just saw a particular way of writing each ordinal $\a\leq\om^\eta$ as a sum of ordinals of the form $\om^{\eta[m_0]...[m_l]}$.
%We will use this decomposition to define $\nabla^{(1+\a)}$ by composing operators of the form $\J^{\om^{\eta[m_0]...[m_l]}}$.
We will now recast this decomposition in terms of composition of operators of the form $\J^{\om^{\eta[m_0]...[m_l]}}$.
%The reason we use this particular decomposition, and not the more standard Cantor normal form, is that otherwise we would not get Lemmas \ref{le: finitely many J} and \ref{le: nabla nested J} below.

\begin{definition} \label{def: nabla}
Given $\eta$ and  $\la n_0,....,n_k\ra\in T_\eta$, we define $\J^{\om^\eta}_{\la n_0,....,n_k\ra}\colon\Bai\to\Bai$ as follows: let $\J^{\om^\eta}_{\la\ra}=\J^{\om^\eta}$ and
\begin{eqnarray*}
\J^{\om^\eta}_{\la n_0,....,n_k\ra} &=&   \J^{\om^{\eta[n_0]}}_{\la n_1,....,n_k\ra} \circ \J^{\om^\eta}_{n_0}  \\
			&=&  \J^{\om^{\eta[n_0]...[n_{k}]}} \circ \bigodot_{i=0}^{n_k-1}\J^{\om^{\eta[n_0]...[n_{k-1}][i]}} \circ\cdots\circ  \bigodot_{i=0}^{n_0-1}\J^{\om^{\eta[i]}},
\end{eqnarray*}
where $\bigodot_{i=0}^l f_i= f^l\circ f^{l-1}\circ\cdots\circ f_0$ is the iterated composition operator.
(Compare with the definition of $\eta\la n_0,...,n_k\ra$ above.)
\end{definition}

It is not hard to see that $\J^{\om^\eta}_{\la n_0,....,n_k\ra}(Z) \teq Z^{(\eta\la n_0,....,n_k\ra)}$.
We will now see that, despite having a complicated definition, the operators $\J^{\om^\eta}_{\la n_0,....,n_k\ra}$ are the ones we want to use as canonical $\eta\la n_0,....,n_k\ra$-th jump operators.
We now define the finite approximations exactly in the same way.

\begin{definition}
Given $\la n_0,....,n_k\ra\in T_\eta$, we define $J^{\om^\eta}_{\la n_0,....,n_k\ra}\colon\Seq\to\Seq$ as follows
\begin{eqnarray*}
J^{\om^\eta}_{\la n_0,....,n_k\ra} &=&   J^{\om^{\eta[n_0]}}_{\la n_1,....,n_k\ra} \circ J^{\om^\eta}_{n_0}  
\end{eqnarray*}
and let $J^{\om^\eta}_{\la\ra}=J^{\om^\eta}$.
\end{definition}

%
%\begin{lemma}\label{le: finitely many J}
%For each $\si\in\Seq$ and $\la n_0,...,n_k\ra\in T_\eta$, $|J^{\om^\eta}_{\la n_0,...,n_k\ra}| < |\si| \dmi (n_0+...+n_k)$.
%\end{lemma}
%\begin{proof}
%We know that for each $\a$ and $\si$, $|J^{\om^\a}(\si)| \leq |\si| \dmi 1$ \cite[Lemma 6.2 (P$^{\om^\a}$5)]{MMVeblen}.
%There are $n_0+\cdots+n_k+1$ applications of functions of the form $J^{\om^\a}$ in the definition of $J^{\om^\eta}_{\la n_0,...,n_k\ra}$.
%The lemma follows by induction.
%\end{proof}

\begin{lemma} \label{le: nabla nested J}
For $\si,\tau\in\Seq$ and $\la m_0,...,m_l\ra<_{KB}\la n_0,...,n_k\ra\in T_\eta$, if $\la\ra\neq J^{\om^\eta}_{\la n_0,...,n_k\ra}(\si)\subseteq J^{\om^\eta}_{\la n_0,...,n_k\ra}(\tau)$, then $J^{\om^\eta}_{\la m_0,...,m_l\ra}(\si)\subseteq J^{\om^\eta}_{\la m_0,...,m_l\ra}(\tau)$.
\end{lemma}
\begin{proof}
Let $h$ be the greatest such that for all $j<h$, $n_j=m_j$ and $n_h\geq m_h$.
So, either $h=k<l$ and $\la n_0,...,n_k\ra= \la m_0,...,m_h\ra$, or $h<k$ and $n_h>m_h$.
Consider the following chain of inequalities. 
\[
\eta\la m_0,...,m_l\ra < \cdots< \eta\la m_0,...,m_h,m_{h+1}\ra < \eta\la m_0,...,m_h\ra \leq \eta\la n_0,...,n_k\ra.
\]
Assuming that  $\la\ra\neq J^{\om^\eta}_{\la m_0,...,m_l\ra}(\si)$ and that 
\begin{align}
\la\ra\neq J^{\om^\eta}_{\la n_0,...,n_k\ra}(\si) &\subseteq J^{\om^\eta}_{\la n_0,...,n_k\ra}(\tau),
\end{align} 
we will show each of the following inclusions one-by-one:
\begin{align}
\la\ra\neq J^{\om^\eta}_{\la m_0,...,m_h\ra}(\si) &\subseteq J^{\om^\eta}_{\la m_0,...,m_h\ra}(\tau)   \label{eq: h} \\
\la\ra\neq J^{\om^\eta}_{\la m_0,...,m_h,m_{h+1}\ra}(\si) &\subseteq J^{\om^\eta}_{\la m_0,...,m_h,m_{h+1}\ra}(\tau) \label{eq: h+1}  \\
 &\vdots&  \nonumber \\ 
\la\ra\neq J^{\om^\eta}_{\la m_0,...,m_h,...,m_l\ra}(\si) &\subseteq J^{\om^\eta}_{\la m_0,...,m_h,...,m_l\ra}(\tau)  \label{eq: l} 
\end{align} 

To prove (\ref{eq: h}) we note, just by unfolding the definitions, that $J^{\om^\eta}_{\la n_0,...,n_k\ra}$ can be written by a composition of a list of operators of the form $J^{\om^\b}$ and $J^{\om^\eta}_{\la m_0,...,m_h\ra}$.
More precisely (although not relevant for the proof), we have
\[
J^{\om^\eta}_{\la n_0,...,n_k\ra} = J^{\om^{\eta[n_0]...[n_h]}}_{\la n_h+1,...,n_k\ra} \circ  \bigodot_{i=m_h+1}^{n_k-1} J^{\om^{\eta[n_k]...[n_{h-1}][i]}} \circ J^{\om^\eta}_{\la m_0,...,m_h\ra}.
\]
An iteration of Corollary \ref{co: omega alpha inclusion} gives us that $J^{\om^\eta}_{\la m_0,...,m_h\ra}(\si) \subseteq J^{\om^\eta}_{\la m_0,...,m_h\ra}(\tau)$.
And since $\la\ra\neq J^{\om^\eta}_{\la n_0,...,n_k\ra}(\si)$, we have that $\la\ra\neq J^{\om^\eta}_{\la m_0,...,m_h\ra}(\si)$ too \cite[Lemma 6.2 (P$^{\om^\a}$1)]{MMVeblen}, proving  (\ref{eq: h}).

Let us now prove (\ref{eq: h+1}).
Let $\basi$, $\batau$ and $\a$ be such that
\[
J^{\om^\eta}_{\la m_0,...,m_h\ra}(\si)= J^{\om^{\a}} (\basi)
	\quad\and\quad
J^{\om^\eta}_{\la m_0,...,m_h\ra}(\tau)= J^{\om^{\a}} (\batau).
\]
More precisely, let $\overline{J^{\om^\eta}_{\la m_0,...,m_h\ra}}$ be defined like $J^{\om^\eta}_{\la m_0,...,m_j\ra}$ but without the application of the first term $ J^{\om^{\eta[m_0]...[m_h]}}$ in the compostion.
That is $\overline{J^{\om^\eta}_{\la m_0,...,m_h\ra}}=\bigodot_{i=0}^{m_h-1}\J^{\om^{\eta[m_0]...[m_{h-1}][i]}} \circ\cdots\circ  \bigodot_{i=0}^{m_0-1}\J^{\om^{\eta[i]}}$,
Let $\basi=\overline{J^{\om^\eta}_{\la m_0,...,m_h\ra}}(\si)$, $\batau = \overline{J^{\om^\eta}_{\la m_0,...,m_h\ra}}(\tau)$, and $\a=\eta[m_0]...[m_h]$.
%
%\nabla^{1+\g}_s=\bigodot_{i=0}^{m_l} J^{\om^{\eta[n_0]...[m_{l-1}][i]}} \circ \cdots\circ  \bigodot_{i=0}^{m_{h+1}} J^{\om^{\eta[m_0]...[m_{h}][i]}}(\nabla^-_s),
%%\nabla^{1+\g}_s=\bigodot_{i=0}^{m_l}J^{\om^{\a[i]}}(\nabla^-_s)  = J^{\om^\a}_{m_l+1}(\nabla^-_s),
%\]

We now have that 
\[
J^{\om^\eta}_{\la m_0,...,m_h,m_{h+1}\ra} (\si)= J^{\om^{\a}}_{m_{h+1}+1} (\basi),
\]
and the same for $\tau$ instead of $\si$.
From (\ref{eq: h}) we have that $\la\ra\neq J^{\om^\a}(\basi) \subseteq J^{\om^\a}(\batau)$.
Then, by Lemma \ref{le: omega alpha inclusion}, $J^{\om^\a}_{m_{h+1}+1}(\basi) \subseteq J^{\om^\a}_{m_{h+1}+1}(\batau)$ as needed for (\ref{eq: h+1}).

Let $\basi_1$, $\batau_1$ and $\a_1$ be such that
\[
J^{\om^\eta}_{\la m_0,...,m_h,m_{h+1}\ra}(\si)= J^{\om^{\a_1}} (\basi_1)
	\quad\and\quad
J^{\om^\eta}_{\la m_0,...,m_h,m_{h+1}\ra}(\tau)= J^{\om^{\a_1}} (\batau_1).
\]
More precisely, $\a_1=\a[m_{h+1}]$, $\basi_1=J^{\om^\a}_{m_{h+1}}(\basi)$ and $\batau_1=J^{\om^\a}_{m_{h+1}}(\batau)$.
Therefore, 
\[
J^{\om^\eta}_{\la m_0,...,m_h,...,m_l\ra} (\si)= J^{\om^{\a}}_{\la m_h+1,...,m_l\ra} (\basi_1).
\]
Since $J^{\om^\eta}_{\la m_0,...,m_h,...,m_l\ra}(\si)\neq\la\ra$, we get that $|\basi_1| \geq 2$ \cite[Lemma 6.2 (P$^{\om^\a}$1)]{MMVeblen}, and hence $J^{\om^\eta}_{\la m_0,...,m_h,m_{h+1}\ra}(\si) = J^{\om^{\a_1}} (\basi_1) \neq\la\ra$ too.
This finishes the proof of (\ref{eq: h+1}).

Continue like this proving the following lines up to (\ref{eq: l}).
\end{proof}

%%%%%%%%%%%%%%%%%%%%%%%%%%%%%%%%%%%%%%%%%%%%%%%%%%%%%%%%%%%%%%%%%%%%%%%%%%%%%%%%%%%%%%%%%%%%%%%%%%%%%%%%%%%%%%%%%%%%%%%%%%%%%%%%%%%%%%%%%%%%%%%%%%%%%%%%%%%%%%%%%%%%%%%%%%%%%%%%%%%%%%%%%%%%%%%%%%%%%%%%%%%%%%%%%%%%%%%%%%%%%%%%%
\section{The pre-orderings of belief}\label{se: belief}

In this section we define the orderings $\leq_\xi$ that we use to define the notion of ``$\xi$-belief,'' or ``looking like a $\xi$-true stage.''

%%%%%%%%%%%%%%%
\subsection{Canonical $\Delta^0_\xi$-complete functions} \label{see: xi-true stages}

For the rest of this section, fix $\eta$ a nice computable ordinal as we had before.

\begin{definition} \label{def: nabla}
Let $\nabla^1\in\Bai$ be the function constant equal to zero.
Consider $\a\in\om^\eta+1$.
Suppose that $\a=\eta\la n_0,...,n_k\ra$. Let
\[
\nabla^{1+\a}=\J^{\om^\eta}_{\la n_0,...,n_k\ra}(\nabla^1).
\]
Recall that if $\a$ is infinite, then $1+\a=\a$.
\end{definition}

The reason we start with $\nabla^1$ rather than $\nabla^0$ is the usual problem with the transfinite arithmetic hierarchy that $0^{(\xi)}$ is $\Delta^0_{\xi+1}$-complete for finite $\xi$ and is $\Delta^0_{\xi}$-complete for infinite $\xi$.
This way we do not have to distinguish between these two cases every time: $\nabla^\xi$ is $\Delta^0_\xi$-complete for all $\xi$ uniformly in $\xi$.

We now define the finite approximations.

\begin{definition}
Consider $s\in\om$.
Let $\nabla^1_s\in\Seq$ be the string with $s$ zeros.
Suppose again that $\a=\eta\la n_0,...,n_k\ra$. Let
\[
\nabla^{1+\a}_s=J^{\om^\eta}_{\la n_0,...,n_k\ra}(\nabla^1_s).
\]
\end{definition}

Let us now prove that these approximations satisfy (\ref{pa:N1})-(\ref{pa:N3}).
For the reader only interested in the finite case, these properties follow easily, and most of the complexity in proving the following lemmas appears only in the infinite case.

\begin{lemma}\label{le: seq of xi true}
For each $\xi$, there is a sub-sequence $\{t_i:i\in\om\}$ such that  $\bigcup_{i\in\om}\nabla^\xi_{t_1} =\nabla^\xi$
\end{lemma}
\begin{proof}
This is just an iteration of Lemma \ref{le: inf true stages} or of \cite[Lemma 6.9]{MMVeblen}.
It follows from taking the inverse of $\J^{\om^\eta}_{\la n_0,...,n_k\ra}$ by composing the appropriate functions $J^{-\om^{\eta[n_0]...,[n_h][i]}}$.
\end{proof}

\begin{lemma}\label{le: finitely many}
For each $s$ there are only finitely many $\xi\leq\om^\eta$ with $\nabla^\xi_s\neq\la\ra$.
Furthermore, we can computably find those $\xi$'s uniformly in $s$.
\end{lemma}
\begin{proof}
We know from \cite[Lemma 6.2 (P$^{\om^\a}$5)]{MMVeblen} that  for each $\a$ and $\si\neq\la\ra$, $|J^{\om^\a}(\si)| \leq |\si| - 1$. 
There are $n_0+\cdots+n_k+1$ applications of functions of the form $J^{\om^\a}$ in the definition of $J^{\om^\eta}_{\la n_0,...,n_k\ra}$.
Thus, we know that for each $\xi=\eta\la n_0,...,n_k\ra$ and $s$, if $s\leq n_0+\cdots+n_k$, then $\nabla^{\xi}_s=\la\ra$.
Now, we just need to notice that for each $s$ there are finitely many $\la n_0,...,n_k\ra \in T_\eta$ with $n_0+\cdots+n_k< s$.
This is because $T_\eta$ is well-founded and for each $i$, $n_i< s$.

Since $T_\eta$ is computable, we can easily list those tuples effectively. 
\end{proof}

\begin{lemma} \label{le: nabla nested}
For $s<t$ and $\g<\xi\leq \om^\eta$, if $\la\ra\neq\nabla^\xi_s\subseteq\nabla^\xi_t$, then $\nabla^\g_s\subseteq\nabla^\g_t$.
\end{lemma}
\begin{proof}
This just follows directly from Lemmas \ref{le: order iso} and \ref{le: nabla nested J}.
\end{proof}

\begin{lemma}\label{le: nabla continuous}
Let $\lambda\leq \om^\eta$ be a limit ordinal and $s<t\in\om$.
Then $\nabla^\lambda_s\subseteq\nabla^\lambda_t$ if and only if $(\forall \xi<\lambda) \nabla^\xi_s\subseteq\nabla^\xi_t$, provided $\la\ra\neq \nabla^\lambda_s$.
\end{lemma}
\begin{proof}
The left-to-right direction follows from the previous lemma.
For the right-to-left direction let $\lambda=\eta\la m_0,...,m_h\ra$.
Let $\basi=\overline{J^{\om^\eta}_{\la m_0,...,m_h\ra}}(\nabla^1_s)$, $\batau = \overline{J^{\om^\eta}_{\la m_0,...,m_h\ra}}(\nabla^1_t)$, and $\a=\eta[m_0]...[m_h]$ (as in the proof of Lemma \ref{le: nabla nested J}).
What matters is that $\basi$, $\batau$ and $\a$ were chosen so that
\[
\nabla^\lambda_s= J^{\om^{\a}} (\basi)
	\quad\and\quad
\nabla^{\eta\la m_0,...,m_h,i\ra}_s= J^{\om^{\a}}_{i+1} (\basi),
\]
and the same for $\batau$ instead of $\basi$.
Since $\eta\la m_0,...,m_h,i\ra<\lambda$, we have that for all $i$, $J^{\om^{\a}}_{i+1} (\basi)\subseteq J^{\om^{\a}}_{i+1} (\batau)$.
It then follows from Lemma \ref{le: omega alpha inclusion} that $J^{\om^{\a}} (\basi)\subseteq J^{\om^{\a}} (\batau)$ as needed.
\end{proof}

%%%%%%%%%%%%%%%%
\subsection{$\xi$-true stages} \label{see: xi-true stages}

The first attempt to define the relation of ``$\xi$-believes in'' is the following.
Given $s,t\in\om$ we define
\[
s\Nleq_\xi t \iff (\forall \g\leq\xi+1)\ \nabla^{\g}_s\subseteq \nabla^{\g}_t.
\]

This definition has almost all the properties we need, except for continuity, which we will fix later in Lemma \ref{le: xi-belief final}.
For finite $\xi$, these relations have all the properties we need, as continuity is not important in that case, and hence there is no need for further modifications.

That each $\Nleq_\xi$ is a pre-ordering is immediate from the definition.
That the sequence is nested is also immediate.
Let us observe that the sequence $(\Nleq_\xi)_{\xi<\om^\eta}$ is computable:
Given $s$ and $t$, by Lemma \ref{le: finitely many}, there are only finitely many $\g$ with $\nabla^\g_s\neq\la\ra$, and we can effectively find them.
We only need to check if  $\nabla^{\g}_s\subseteq \nabla^{\g}_t$  among those $\g$'s.
We also observe that  $\Nleq_0=\leq$, and that for each $\a$, the $t$'s for which $\nabla^{\a+1}_t$ is correct and non-empty form a $\Nleq_\a$-ascending sequence.

\begin{lemma} \label {le: Nleq}
$(\Nleq_\xi)_{\xi<\om^\eta}$ is a nested computable sequence of pre-orderings satisfying $(\clubsuit)$.
\end{lemma}
\begin{proof}
We already pointed out that it is a nested computable sequence of pre-orderings.
To prove $(\clubsuit)$, consider $r<s<t$ with $r\Nleq_{\xi+1}t$ and $s\Nleq_\xi t$--we want to show that $r\Nleq_{\xi+1}s$.
Suppose, toward a contradiction that $r\not\Nleq_{\xi+1}s$.
Then, for some $\g\leq \xi+2$ we have that $\nabla^{\g}_r\not\subseteq\nabla^\g_s$.
We may assume $\g$ is a successor, as if $\nabla^{\lambda}_r\not\subseteq\nabla^\lambda_s$ for a limit ordinal $\lambda$, then $\nabla^{\g}_r\not\subseteq\nabla^\g_s$ for some $\g<\lambda$ (as it follows from Lemma \ref{le: nabla continuous}).
Let $\d=\g-1$.
So we have that $J(\nabla^\d_r)\not\subseteq J(\nabla^\d_s)$ and, since $\d\leq\xi$, that $\nabla^\d_r\subseteq\nabla^\d_t$ and $\nabla^\d_s\subseteq\nabla^\d_t$.
Therefore $\nabla^\d_r\subseteq\nabla^\d_s\subseteq\nabla^\d_t$.
From Lemma \ref{le: s t u} we get that $J(\nabla^\d_r)\not\subseteq J(\nabla^\d_t)$ and hence that $r\not\Nleq_{\xi+1}t$ giving the desired contradiction. 
\end{proof}

Now we will modify these relations slightly to get continuity.
Just to describe the intuition behind these modifications, let us assume $\eta=\om$.
The problem we have is that we might have $\nabla^{\om+1}_s\not\subseteq\nabla^{\om+1}_t$, and hence $s\not\Nleq_\om t$ while  $\nabla^{\om}_s\subseteq\nabla^{\om}_t$, which is equivalent to $(\forall n<\om)\ \nabla^{n}_s\subseteq\nabla^{n}_t$.
Let's look at this more carefully.
Let $j$ be the greatest such that $\nabla^{\om+1}_s\upto j\subseteq\nabla^{\om+1}_t$.
So, we must have that $\varphi_j^{\nabla^{\om}_s}(j)\diverges$ but $\varphi_j^{\nabla^{\om}_t}(j)\converges$.
But we should not blame $\om$ for this disagreement.
Let $t_j$ be as in the definition of $\nabla^{\om}_t$-true stages (Definition \ref{def:1true}, so $t_j$ is essentially the oracle-use of the computation $\varphi_j^{\nabla^{\om}_t}(j)\converges$), so that $\nabla^{\om+1}_t(j)=\ulcorner \nabla^{\om}_t\upto t_j\urcorner = \ulcorner \la \nabla^{1}_t(0),\nabla^2_t(0),...\nabla^{t_j}_t(0)\ra\urcorner$.
The one we should blame is $t_j$, as it is the one witnessing a computation that we thought should have diverged, and thus, we should let $s\not\leq_{t_j} t$.

\begin{lemma}\label{le: xi-belief final}
There is a sequence  $(\leq_\xi)_{\xi<\om^\eta}$ satisfying (\ref{pa: 0})-(\ref{pa:6}) and $(\clubsuit)$.
\end{lemma}
\begin{proof}
We will modify the relations $\Nleq_\xi$ to get them to satisfy continuity.
Let us start by defining  a list of all the possible situations that cause a problem for continuity.
Let $C$ be the set of all triples $(\lambda,u,v)$ where $\lambda< \om^\eta$ is a limit ordinal, $\nabla^{\lambda}_u\subsetneq\nabla^{\lambda}_v$, $\nabla^{\lambda+1}_u\not\subseteq\nabla^{\lambda+1}_v$, and  $v$ is minimal in the sense that if $\nabla^{\lambda}_u\subsetneq\nabla^{\lambda}_r\subsetneq\nabla^{\lambda}_v$, then $\nabla^{\lambda+1}_u\subseteq\nabla^{\lambda+1}_r$.
For each such tuple, let $\g_{\lambda,v}$ be such that the last entry of $\nabla^{\lambda}_v$ is $\nabla^{\g_{\lambda,v}}_v(0)$, (or, to be more specific $\g_{\lambda,v}=1+\eta\la n_0,...,n_k,j-1\ra$ where $\lambda=\eta\la n_0,...,n_k\ra$ and $j=|\nabla^{\lambda}_v|$).
(Recall that the idea is to blame $\g_{\lambda,v}$ for the problem rather than $\lambda$, and let $u\not\leq_{\g_{\lambda,v}+1}v$.)
Define
\[
s\leq_{\xi} t  \quad\iff 	\quad  s\Nleq_\xi t \ \ \ \& \ \ \ \neg \exists (\lambda,u,v)\in C \left(\g_{\lambda,v}< \xi \  \ \&\ \  u \leq s < v \Nleq_{\g_{\lambda,v}} t\right).
\]
Let us now show that this sequence satisfies all the desired properties:

Is easy to see that (\ref{pa: 0}) holds, i.e., that $\leq_0=\leq$ . 

{\em Computability (\ref{pa: 0.5}).}
It is computable because $\Nleq$ is, and because the existential quantifier over $(\lambda,u,v)\in C$ is bounded, as $u,v\leq t$, and for each $v$, there are only finitely many $\lambda$'s with $\nabla^{\lambda+1}_v\neq\la\ra$.

{\em Pre-ordered (\ref{pa: 1}).} Reflexivity follows from the reflexivity of $\Nleq$, and the fact that if $s=t$, then there is no $v$ with $s<v\leq t$.
For transitivity suppose that $s\leq_\xi t\leq_\xi r$ but that $s\not\leq_\xi r$.
By the transitivity of $\Nleq$ we do have that $s\Nleq_\xi r$, so there must be a $(\lambda,u,v)\in C$ witnessing that $s\not\leq_\xi r$.
If $v>t$, then $u\leq t<v\Nleq_{\g_{\lambda,v}}r$, so $(\lambda,u,v)$ also witnesses that $t\not\leq_\xi r$ giving a contradiction.
If $v\leq t$, then since $v\Nleq_\xi r$ and $t\Nleq_\xi r$, by $(\clubsuit)$ applied to $\Nleq$, we get that $v\Nleq_\xi t$ (using Observation \ref{rm: tree like}).
Since $\g_{\lambda,v}<\xi$, we have that  $u\leq s<v\Nleq_{\g_{\lambda,v}}t$, and hence $(\lambda,u,v)$ witnesses that $s\not\leq_\xi t$ again giving a contradiction.

{\em Nested (\ref{pa: 2}). } Showing that it is nested is quite straightforward.

{\em Continuity (\ref{pa:4}).}
Suppose, towards a contradiction, then for some limit $\a$ we have $s\not\leq_\a t$ but that $(\forall \xi<\a)\ s\leq_\xi t$.
If $s\not\leq_\a t$ due to some $(\lambda, u,v)\in C$, then we would also have $s\not\leq_{\g_{\lambda,v}+1}t$ and $\g_{\lambda,v}<\a$.
So it must be that $s\not\Nleq_\a t$.
Since $(\forall \xi<\a)\ s\Nleq_\xi t$, we have that $\nabla^\lambda_s\subseteq\nabla^\lambda_t$ and $\nabla^{\lambda+1}_s\not\subseteq\nabla^{\lambda+1}_t$.
Let $v$ be the least such that $\nabla^\lambda_s\subsetneq\nabla^\lambda_v\subseteq\nabla^\lambda_t$ and $\nabla^{\lambda+1}_s\not\subseteq\nabla^{\lambda+1}_v$.
We then have that $(\lambda,s,v)\in C$.
We also have that $v\Nleq_{\g_{\lambda,v}} t$ because $\la\ra\neq\nabla^\lambda_v\subseteq\nabla^\lambda_t$.
Therefore, $s\not\leq_{\g_{\lambda,v}+1} t$ contradicting our assumptions.

It is easy to see that (\ref{pa:5}) holds, i.e., that $s\leq_\xi t \implies  \nabla^{\xi+1}_s\subseteq\nabla^{\xi+1}_t$ for all $\xi$.

Suppose  that $\la\ra\neq \nabla^{\om^\eta}_s\subseteq\nabla^{\om^\eta}_t$--we want to show that $s\leq_\xi t$ for all $\xi<\eta$.
Part (\ref{pa:6}) would then follow from Lemma \ref{le: seq of xi true}.
From Lemma  \ref{le: nabla nested} we  have that $s\Nleq_\xi t$ for all $\xi<\om^\eta$.
Take $(\lambda,u,v)\in C$ and suppose towards a contradiction that $u \leq s < v \Nleq_{\g_{\lambda,v}} t$.
Since $v \Nleq_{\g_{\lambda,v}} t$,  $\nabla^{\g_{\lambda,v}}_v(0)=\nabla^{\g_{\lambda,v}}_t(0)$, and since $\nabla^{\g_{\lambda,v}}_v(0)$ is the last entry of $\nabla^{\lambda}_v$, using Lemma \ref{le: seq of xi true} we have that $\nabla^{\lambda}_v\subseteq\nabla^{\lambda}_t$ because all the entries of $\nabla^{\lambda}_v$ are then the same.
Since we are assuming that $\la\ra\neq \nabla^{\om^\eta}_s\subseteq\nabla^{\om^\eta}_t$, we also have that $\nabla^{\lambda}_s\subseteq\nabla^{\lambda}_t$.
From the fact that $(\lambda,u,v)\in C$ we have that $\nabla^{\lambda}_v\subseteq\nabla^{\lambda}_v$.
So, given that all these strings are compatible we  have that $\nabla^{\lambda}_u\subseteq\nabla^{\lambda}_s\subseteq\nabla^{\lambda}_v\subseteq\nabla^{\lambda}_t$.
By the minimality of $v$ in the definition of $C$ we have that $\nabla^{\lambda+1}_u\subseteq\nabla^{\lambda+1}_s\not\subseteq\nabla^{\lambda+1}_v$.
From Lemma \ref{le: s t u} we get that $\nabla^{\lambda+1}_s\not\subseteq\nabla^{\lambda+1}_t$, and hence that $s\not\Nleq_{\eta} t$ giving the desired contradiction.

$(\clubsuit)$. Suppose that $s<t<r$, $s\leq_{\xi+1} r$ and $t\leq_\xi r$ but that $s\not\leq_{\xi+1} t$.
Since $(\clubsuit)$ holds for $(\Nleq_\xi)_{\xi<\om^\eta}$, we have that $s\Nleq_{\xi+1} t$, so there must be a  $(\lambda,u,v)\in C$ witnessing  that $s\not\leq_{\xi+1} t$.
Thus $v\Nleq_{\g_{\lambda,v}} t \Nleq_\xi r$.
Since $\xi+1>{\g_{\lambda,v}}$ we get that $v\Nleq_{\g_{\lambda,v}} r$.
So $(\lambda,u,v)$ witnesses that $s\not\leq_{\xi+1} r$ giving a contradiction.
\end{proof}

\bibliography{bftypes}
\bibliographystyle{alpha}

\end{document}